# Families of Linearly χ-bounded Graphs without Chair or its Induced Subgraphs

**Medha Dhurandhar**

**Abstract**: A hereditary class $\mathscr{H}$ of graphs is χ-bounded if there is a χ-binding function f s.t. for every G ∈ $\mathscr{H}$, $\chi(G) \leq f(\omega(G))$. Here we prove that if a graph G is free of 1. {Chair; $P_4+K_1$} or 2. {Chair; HVN}, then χ(G) is linearly bounded by maximum clique size of G. We further prove that if G is free of 3. {$P_4+K_1$; $P_3\cup K_1$} or 4. {$P_4+K_1$; $K_2\cup 2K_1$} or 5. {HVN; $P_3\cup K_1$} or 6. {HVN; $K_2\cup 2K_1$} or 7. {$K_5$-e; $P_3\cup K_1$} or 8. {$K_5$-e; $K_2\cup 2K_1$}, then there is a tight linear χ-bound for G.

**Introduction:**
Graph coloring is one of the best known, popular and extensively researched subject in the field of graph theory, having many applications and conjectures, which are still open and studied by various mathematicians and computer scientists along the world. In general, there is no bound on the chromatic number of a graph in terms of the size of its largest complete subgraph, since there are graphs containing no triangle, but having arbitrarily large chromatic number. In [2], [3], [4], [5], [6] we find surveys on vertex chromatic number of a graph.

We consider here finite, simple and undirected graphs. For terms which are not defined herein we refer to Bondy and Murty [1].

**Notation:** For a graph G, V(G), E(G), Δ(G), ω(G), χ(G) denote the vertex set, edge set, maximum degree, size of a maximum clique, chromatic number respectively. For u ∈ V(G), N(u) = {v ∈ V(G) /uv ∈ E(G)} and $\overline{N(u)}$ = N(u)∪(u). If S ∈ V(G), then <S> denotes the subgraph of G induced by S. We refer to [] for an extensive survey of χ-bounds for various classes of graphs. Here we are interested in classes of graphs that admit a linear χ-binding function.

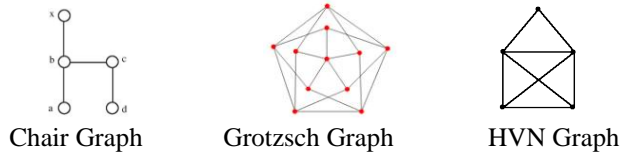

Chair Graph    Grotzsch Graph    HVN Graph
**Figure 1**

We state below a **Lemma,** which is used in all following Theorems.

**Lemma**: Let G be a χ-critical, Chair-free graph. Among all (χ-1)-colorings of G-v for every v ∈ V(G), let ∃ u ∈ V(G) and (χ-1)-coloring C of G-u s.t. u has maximum no. |R| vertices with unique colors and if u has repeat colors $\alpha_1, ..., \alpha_{\chi-|R|-1}$, then u has minimum no. $N_1$ of $\alpha_1$-vertices. Further with |R| unique colors and $N_1$ no. of $\alpha_1$-vertices, let u have minimum no. $N_2$ of $\alpha_2$-vertices and so on. Let R = {x ∈ N(u)/ x receives a unique color from C in <N(u)>}. Then
**A.** Every x ∈ R has a $\alpha_i$-vertex in <N(u)> ∀ i = 1, 2,.., χ-|R|-1.
**B.** Every $\alpha_{ij}$ ∈ N(u)-R has a $\alpha_k$-vertex in <N(u)> for k > i.
Proof:
   **A.** Let if possible ∃ x ∈ R s.t. x has no $\alpha_i$-vertex in <N(u)> for some i ∈ {1, 2,.., χ-|R|-1}. If x has $\alpha_i$-vertex x' ∈ <V(G)-N(u)>, then <u, $\alpha_{i1}, \alpha_{i2}$, x'> = Chair and if x has no $\alpha_i$-vertex, then color by x by $\alpha_i$, u by color of x, a contradiction in both the cases.
   **B.** Let if possible ∃ $\alpha_{ij}$ ∈ N(u)-R s.t. $\alpha_{ij}$ has no $\alpha_k$-vertex in <N(u)> for some i, j, k ∈ {1, 2,.., χ-|R|-1} with k > i. As G is Chair-free as in **A**, $\alpha_{ij}$ has no $\alpha_k$-vertex in G. Now color $\alpha_{ij}$ by $\alpha_k$ and let C' be this (χ(G)-1)-coloring of G-u. Now $N_i > 2$ (else in C', u has |R|+1 unique colors) ⇒ <N(u)>

uses |R| unique colors in C', $N_s$ no. of $\alpha_s$ colors for s < i and $N_i$-1 no. of $\alpha_i$ colors, a contradiction to our assumption.

This proves the **Lemma**.

**Theorem 1:** If G is {Chair, $(P_4+K_1)$}-free, then $\chi(G) \leq 2\omega-1$.
Proof: Let if possible G be a smallest {Chair, $(P_4+K_1)$}-free graph with $\chi(G) > 2\omega-1 \Rightarrow \chi(G)-1 \leq \chi(G-u) \leq 2\omega-1 \,\forall\, u \in V(G)$. Among all ($\chi$-1)-colorings of G-v for every v ∈ V(G), let ∃ u ∈ V(G) and ($\chi$-1)-coloring C of G-u s.t. u has maximum no. |R| vertices with unique colors and if u has repeat colors $\alpha_1, ..., \alpha_{\chi-|R|-1}$, then u has minimum no. $N_1$ of $\alpha_1$-vertices. Further with |R| unique colors and $N_1$ no. of $\alpha_1$-vertices, let u have minimum no. $N_2$ of $\alpha_2$-vertices and so on.

Let R = {x ∈ N(u)/ x receives a unique color from C in <N(u)>}. Let <$R_1$> be a maximum clique in <R>. By **Lemma** every x ∈ $R_1$ has a $\alpha_i$-vertex in <N(u)> and as G is $(P_4+K_1)$-free, clearly ∃ $T_1 = \{\alpha_{1i}, \alpha_{2k}, ..., \alpha_{(\chi-|R|-1)j}\} \subseteq N(u)-R$ s.t. xy ∈ E(G) ∀ x ∈ $R_1$ and y ∈ $T_1$. If $\alpha_{mi}\alpha_{nk} \notin E(G)$ for some m < n where $\alpha_{mi}, \alpha_{nk} \in T_1$, then by **Lemma** ∃ l s.t. $\alpha_{mi}\alpha_{nl} \in E(G)$. Now $x\alpha_{nl} \in E(G)$ ∀ x ∈ $R_1$ (else if $x\alpha_{nl} \notin E(G)$, then <$\alpha_{nl}, \alpha_{mi}, x, \alpha_{nk}, u$> = $(P_4+K_1)$). Let $T_2 = \{T_1-\alpha_{nk}\}\cup\alpha_{nl}$ and so on. Thus ∃ $T_m$ s.t. <$R_1\cup T_m\cup\{u\}$> is complete and $|R_1|+\chi-|R|-1+1 \leq \omega \Rightarrow \chi \leq \omega+|R|-|R_1|$. **I**

**Claim**: <R-$R_1$> is complete.
Let if possible ∃ x, y ∈ R-$R_1$ s.t. xy ∉ E(G). As G is $(P_4+K_1)$-free, ∃ a ∈ $R_1$ s.t. ax, ay ∉ E(G) (else let ax ∉ E(G) and ay ∈ E(G). Then ∃ b ∈ $R_1$ s.t. by ∉ E(G) ⇒ bx ∈ E(G) and <x, b, a, y, u> = $P_4+K_1$). Let ax', ay', xa' ∈ E(G). As G is Chair-free, x'y, y'x, a'y ∈ E(G) and y' (a') is a unique vertex of a, x (y, x) of that color. Also y'a' ∈ E(G) (else y''a'' ∈ E(G) and <y', a, a'', x, a'> = Chair). Similarly, x'y, x'y', x'a' ∈ E(G) and <a, y', a', y, x'> = $(P_4+K_1)$, a contradiction. Hence the **Claim** holds.

⇒ |R-$R_1$| ≤ $\omega$-1 and by **I**, $\chi \leq 2\omega-1$, a contradiction.

This proves **Theorem 1**.

**Corollary 1**: If G is {$(P_4+K_1)$, $K_{1,3}$}-free, then $\chi(G) \leq 2\omega-1$.

**Corollary 2**: If G is {$(P_4+K_1)$, $(P_3\cup K_1)$}-free, then $\chi(G) \leq \frac{3\omega}{2}$.

Proof: Let R = {x ∈ N(u)/ x receives a unique color from C in <N(u)>}. Let <$R_1$> be a maximum clique in <R>. By **Lemma** every x ∈ $R_1$ has a $\alpha_i$-vertex in <N(u)> and as G is $(P_4+K_1)$-free, clearly ∃ $T_1 = \{\alpha_{1i}, \alpha_{2k}, ..., \alpha_{(\chi-|R|-1)j}\} \subseteq N(u)-R$ s.t. xy ∈ E(G) ∀ x ∈ $R_1$ and y ∈ $T_1$. If $\alpha_{mi}\alpha_{nk} \notin E(G)$ for some m < n where $\alpha_{mi}, \alpha_{nk} \in T$, then ∃ l s.t. $\alpha_{mi}\alpha_{nl} \in E(G)$. Now $x\alpha_{nl} \in E(G)$ ∀ x ∈ $R_1$ (else if $x\alpha_{nl} \notin E(G)$, then <u, $\alpha_{nl}, \alpha_{mi}, x, \alpha_{nk}$> = $(P_4+K_1)$). Let $T_2 = \{T_1-\alpha_{nk}\}\cup\alpha_{nl}$ and so on. Thus ∃ $T_m$ s.t. <$R_1\cup T_m\cup\{u\}$> is complete and $|R_1|+\chi-|R|-1+1 \leq \omega \Rightarrow \chi \leq \omega+|R|-|R_1|$. **I**

Let S = {x' ∈ V(G)-N(u)/ ∃ x ∈ R with same color s.t. xa ∉ E(G) for some a ∈ R}. Since $(P_3\cup K_1)$ is an induced subgraph of Chair, from **Theorem 1**, it follows that <R-$R_1$> is complete.

**Claim 1**: <S> is complete.
Let if possible ∃ v', w' ∈ S s.t. v'w' ∉ E(G). By construction of S, ∃ v ∈ R with same color as v' and a ∈ R with av ∉ E(G) ⇒ av' ∈ E(G) (else <v, u, a, v'> = $P_3\cup K_1$) ⇒ aw' ∈ E(G) (else <v', a, u, w'> = $P_3\cup K_1$) and a ≠ w ⇒ vw' ∈ E(G) (else <v', a, w', v> = $P_3\cup K_1$). But then <u, v, w', v'> = $P_3\cup K_1$, a contradiction. Hence **Claim 1** holds.

**Claim 2**: $|R| \leq \omega$.

Let if possible $|R| > \omega \Rightarrow <R>$ is not complete and $|S| > 1$. Let $T = \{x \in R/ xv \in E(G) \; \forall \; v \in R\}$. Now $|T| > 0$ (else $|S| = |R|$, and by **Claim 1** $|R| \leq \omega$). Let $x \in T$ and $a' \in S \Rightarrow xa' \notin E(G)$ (else let a, b $\in$ S s.t. ab $\notin$ E(G), then $<a, u, b, a', x> = P_4+K_1$). Now a' has a vertex say x' $\in$ V(G)-N(u) of color of x (else color a' by color of x; b by color of a; u by color of b). As G is ($P_3 \cup K_1$)-free, clearly $x'z' \in E(G)$ $\forall \; z' \in S$. Let $T' = \{x' \in V(G)-N(u)/ \exists \; x \in T$ with same color and $x's' \in E(G) \; \forall \; s' \in S\}$. Then as G is ($P_3 \cup K_1$)-free, $<S \cup T'>$ is complete $\Rightarrow |S \cup T'| = |R|$ and **Claim 2** holds.

Now $|R-R_1| \leq |R_1|$, hence $|R-R_1| \leq \frac{\omega}{2}$ and by **I**, $\chi \leq \frac{3\omega}{2}$.

This proves **Corollary 2**.

**Corollary 3**: If G is $\{(P_4+K_1), (K_2 \cup 2K_1)\}$-free, then $\chi(G) \leq \frac{3\omega}{2}$.

Proof: Let $R = \{x \in N(u)/ \; x$ receives a unique color from C in $<N(u)>\}$ and $<R_1>$ be a maximum clique in $<R>$. Now by **Lemma**, every $x \in R_1$ has a $\alpha_i$-vertex in N(u) and as G is $(P_4+K_1)$-free, clearly $\exists \; T_1 = \{\alpha_{1i}, \alpha_{2k}, ..., \alpha_{(\chi-|R|-1)j}\} \subseteq N(u)-R$ s.t. $x\alpha_{in} \in E(G) \; \forall \; x \in R_1$ and $\alpha_{in} \in T_1$. Now let if possible $\alpha_{1i}\alpha_{2k} \notin E(G)$, then by **Lemma**, $\exists$ n s.t. $\alpha_{1i}\alpha_{2n} \in E(G)$. Now $x\alpha_{2n} \in E(G) \; \forall \; x \in R_1$ (else let $x\alpha_{2n} \notin E(G)$ and $<u, \alpha_{2k}, x, \alpha_{1i}, \alpha_{2n}> = (P_4+K_1)$). Let $T_2 = \{T_1 - \alpha_{2k}\} \cup \alpha_{2n}$ and so on $\Rightarrow \exists \; T_m$ s.t. $<R_1 \cup T_m \cup \{u\}> = K_{|R_1|+\chi-|R|-1+1} \Rightarrow \chi \leq \omega+|R-R_1| \Rightarrow |R-R_1| > 1$. **I**

Let $S = \{x' \in V(G)-N(u)/ \exists \; x \in R$ with same color and $ax \notin E(G)$ for some $a \in R\}$. Since $(K_2 \cup 2K_1)$ is an induced subgraph of Chair, from **Theorem 1**, it follows that $<R-R_1>$ is complete.

**Claim 1**: $<S>$ is complete.
Let if possible $\exists \; v', w' \in S$ s.t. $v'w' \notin E(G)$. Let v, w be the corresponding vertices in R $\Rightarrow$ As G is $(K_2 \cup 2K_1)$-free, clearly vw', wv', vw $\in$ E(G). By construction of S $\exists \; z \in R$ s.t. vz $\notin$ E(G). Let z' $\in$ S. Then v'z', vz', zv' $\in$ E(G). Also w'z' $\in$ E(G) (else $<v', z', u, w'> = K_2 \cup 2K_1$) $\Rightarrow$ wz' E(G) (else color v by color of z; z' by color of w; and w', u by color of v). But then $<w, u, v, z', v'> = P_4+K_1$, a contradiction. Hence **Claim 1** holds.

Next let $T = \{x \in R/ xv \in E(G) \; \forall \; v \in R\}$. Let x $\in$ T, a' $\in$ S and a, b $\in$ R s.t. ab $\notin$ E(G). Now a'x $\notin$ E(G) (else $<x, a, u, b, a'> = P_4+K_1$). Similarly, b'x $\notin$ E(G) and a' (b') has a vertex say x' of color of x (else color a' by color of x; b by color of a; u by color of b). As G is $K_2 \cup 2K_1$-free, x' is the only such vertex in V(G)-N(u) $\Rightarrow$ x'b' $\in$ E(G). Thus if $T' = \{x' \in V(G)-N(u) / \exists \; x \in T$ with same color$\}$, then every vertex of T' is adjacent to every vertex of S.

Finally to show that $<T'>$ is complete. Let if possible $\exists \; x', y' \in T'$ s.t. $x'y' \notin E(G)$. As $|R-R_1| > 1$, $\exists$ a, b $\in$ R s.t. ab $\notin$ E(G). As before x, y are non-adjacent to a', b'. Color y by color of a, a by color of b; b' by color of x; x' and u by color of y, a contradiction.

Thus $<T' \cup S> = K_{|R|} \Rightarrow |R| \leq \omega$. **II**

By **I**, $\chi \leq \omega+|R-R_1|$, $|R-R_1| \leq |R_1|$ and by **II**, $|R| \leq \omega \Rightarrow \chi \leq \frac{3\omega}{2}$.

This proves **Corollary 3.**

**Theorem 2:** If G is {Chair, HVN}-free, then $\chi(G) \leq \frac{3\omega}{2}$.

Proof: Let if possible G be a smallest {Chair, HVN}-free graph with $\chi(G) > \frac{3\omega}{2}$. Among all $(\chi-1)$-colorings of G-v for every v $\in$ V(G), let $\exists$ u $\in$ V(G) and $(\chi-1)$-coloring C of G-u s.t. u has maximum no. |R| vertices with unique colors and if u has repeat colors $\alpha_1, ..., \alpha_{\chi-|R|-1}$, then u has minimum no. $N_1$ of $\alpha_1$-vertices. Further with |R| unique colors and $N_1$ no. of $\alpha_1$-vertices, let u have minimum no. $N_2$ of

$\alpha_2$-vertices and so on. Let $R = \{x \in N(u) / $ x receives a unique color from C in $\langle N(u) \rangle\}$. Let $\langle R_1 \rangle$ be a maximum clique in $\langle R \rangle$. Let $S = \{x' \in V(G) - N(u) / \exists\ x \in R$ with same color$\}$.

First we prove that $|R| \leq \omega$. If $|R_1| = 1$, then as G is Chair-free, $\langle S \rangle = K_{|R|} \Rightarrow |R| \leq \omega$. Again if $|R - R_1| = 1$, then $|R| \leq \omega$. Hence let $|R_1| > 1$ and $|R - R_1| > 1$.

**Case 1**: $\exists\ x, y \in |R - R_1|$ s.t. $xy \notin E(G)$.
**Case 1.1**: $\exists\ a \in R_1$ s.t. ax, ay $\notin E(G)$.
Let ay', ya' $\in E(G)$. Then as G is Chair-free y'x, a'x $\in E(G)$ and y', a' are unique vertices of x of those colors. Also y'a' $\in E(G)$ (else $\exists$ a" of color of a s.t. y'a" $\in E(G)$ and $\langle$y', a, a", x, a'$\rangle$ = Chair). Let ax' $\in E(G) \Rightarrow$ yx' $\in E(G)$. Then x' (y', a') is a unique vertex of a and y (a and x; x and y) of that color $\Rightarrow$ x'y', x'a' $\in E(G)$. Let b $\in R_1$. As G is HVN-free, w.l.g. let by' $\notin E(G)$.
**Case 1.1.1**: Let y' have a unique vertex b' of color of b.
Now a'b' $\in E(G)$ (else ab' $\in E(G)$ {else $\langle$y', b', a', a, u$\rangle$ = Chair} $\Rightarrow$ b'x' $\notin E(G)$ {else $\langle$a, b', y', x', a'$\rangle$ = HVN} $\Rightarrow \langle$a, u, b', x', a'$\rangle$ = Chair). Similarly, x'b' $\in E(G) \Rightarrow$ b'a $\in E(G)$ (else $\langle$b', a', x', y', a$\rangle$ = HVN) $\Rightarrow$ Similarly, xb' $\in E(G)$. But then $\langle$a, x', b', y', x$\rangle$ = HVN, a contradiction.
**Case 1.1.2**: Let y'b', y'b" $\in E(G)$ where b', b" $\in V(G) - N(u)$ and have same color as b.
$\Rightarrow$ a (x) is adjacent to either b' or b" (else $\langle$y', b', b", a (x), u$\rangle$ = Chair). Now a, x both are adjacent to either b' or b" (else let ab', xb" $\in E(G) \Rightarrow$ a'b" $\in E(G)$ {else $\langle$y', a', b", a, u$\rangle$ = Chair}. Similarly, x'b" $\in E(G)$ and $\langle$b", a', x', y', a$\rangle$ = HVN). Let ab', xb' $\in E(G) \Rightarrow$ b'x' $\notin E(G)$ (else $\langle$a, x', b', y', x$\rangle$ = HVN). Similarly, b'a' $\notin E(G)$. But then $\langle$a, u, b', x', a'$\rangle$ = Chair, a contradiction.
**Case 1.2**: **Case 1.1** does not exist.
Let ay $\notin E(G) \Rightarrow$ ax $\in E(G)$. Let xb $\notin E(G) \Rightarrow$ yb $\in E(G)$. Now $|R_1| = 2$ (else if $\exists$ c $\in R_1$, then as G is HVN-free xc, yc $\in E(G)$ and $\langle$x, a, c, u, y$\rangle$ = HVN). Let if possible x, a have no common vertex of color of y. Let xy', ay" $\in E(G)$. Now b is not adjacent to both y' and y" (else $\langle$b, y, y", y', x$\rangle$ = Chair). Also by" $\in E(G)$ (else $\langle$a, y", x, b, y$\rangle$ = Chair). Let bx' $\in E(G)$. Now x'y" $\in E(G)$ (else $\langle$b, y", x', u, x$\rangle$ = Chair) $\Rightarrow$ As G is HVN-free, ax' $\notin E(G)$ and x is a unique vertex of a of that color $\Rightarrow$ yx' $\in E(G)$ (else $\langle$b, y, x', a, x$\rangle$ = Chair) and if P is a bicolor y-x path, then $\exists$ x" with same color as x, adjacent to y" on P and $\langle$y", x', x", a, x$\rangle$ = Chair, a contradiction. Hence x, a have a common vertex y' of color of y and y, b have a common vertex x' of color of x. Let xb' $\in E(G)$ where b' has same color as b and ya' $\in E(G)$ where a' has same color as a. Now either b'y' or a'x' $\in E(G)$ (else color x by color of b; b' by color of y; y by color of a; a', u by color of x). W.l.g. let b'y' $\in E(G) \Rightarrow$ As G is HVN-free, ab' $\notin E(G)$ and b is a unique vertex of a of that color $\Rightarrow$ a'x' $\in E(G)$ (else $\exists$ a" ($\neq$ a) of color of a s.t. x'a" $\in E(G)$ {else color a by color of b; b by color of x; x', u by color of a} and $\langle$y, u, a', x', a"$\rangle$ = Chair) $\Rightarrow$ ba' $\notin E(G)$ and a is a unique vertex of b of that color. Now x has a vertex say z ($\neq$ a) of color of a (else color x, b by color of a; a by color of b; u by color of x) $\Rightarrow$ z = a' and xa' $\in E(G)$ (else let xa" $\in E(G)$ where a" has color of a $\Rightarrow$ b'a" $\in E(G)$ {else $\langle$x, a", b', a, b$\rangle$ = Chair} $\Rightarrow$ y'a" $\notin E(G)$ and $\langle$x, a", y', u, y$\rangle$ = Chair) $\Rightarrow$ b'a' $\in E(G)$ (else $\langle$x, b', a', a, b$\rangle$ = Chair) $\Rightarrow$ y'a' $\notin E(G)$ (else $\langle$a', b', x, y', a$\rangle$ = HVN) $\Rightarrow \exists$ a" of color of a s.t. y'a" $\in E(G)$. Now xa" $\notin E(G)$ (else $\langle$x, a", a', a, b$\rangle$ = Chair) $\Rightarrow \langle$x, u, a', y', a"$\rangle$ = Chair, a contradiction.

**Case 2**: $\langle R - R_1 \rangle$ is complete. Let x, y $\in R - R_1$.
W.l.g. bx, ay $\notin E(G)$. Then either yb, xa both $\notin E(G)$ or both $\in E(G)$ (else we get **Case 1.2**).
**Case 2.1**: yb, xa both $\in E(G)$.
Let ay', bx', ya', xb' $\in E(G)$ where a', b', x', y' $\in V(G) - N(u)$ and have colors of a, b, x, y resply.
**Claim**: a has a unique vertex of color of either b or x.
Let if possible a have multiple vertices of color of b and x $\Rightarrow$ ab' $\in E(G)$ (else if ab" $\in E(G)$, then $\langle$a, b", b, x, b'$\rangle$ = Chair). Similarly, ax' $\in E(G)$.
**a)**  by' $\in E(G)$.
$\Rightarrow$ x'y' $\notin E(G)$ (else $\langle$x', y', a, b, u$\rangle$ = HVN) $\Rightarrow$ y'x (x'y) $\in E(G)$ (else $\langle$b, y', x', u, x$\rangle$ = Chair$\rangle$) and clearly, x (y) is a unique vertex of y' (x') of that color $\Rightarrow$ y'b' $\notin E(G)$ and as before y (b) is a unique vertex of b' (y') of that color $\Rightarrow$ ba' $\in E(G)$, a'x' $\notin E(G)$ and x is a unique vertex of a' of that color. But then color x by color of b; b' by color of y; y by color of a; a', u by color of x, a contradiction.

**b)** by', xy' $\notin$ E(G).

$\Rightarrow$ y is a unique vertex of b, x of that color $\Rightarrow$ yb', yx' $\in$ E(G) (else color x, b by color of y; y by color of b; u by color of x). First let a'x $\notin$ E(G) $\Rightarrow$ a is a unique vertex of x of that color (else if xa" $\in$ E(G) where a", a have same color, then <x, a, a", y, a'> = Chair). Now ba' $\notin$ E(G) (else x'a' $\in$ E(G) {else <b, x', a', u, x> = Chair} and <x', a', b, y, u> = HVN) $\Rightarrow$ a is a unique vertex of b of that color and y'a' $\notin$ E(G) (else <a, b, x, y', a'> = Chair) $\Rightarrow$ $\exists$ a" of color of a adjacent to y' on a-y bicolor path $\Rightarrow$ <a, b, x, y', a"> = Chair, a contradiction. Hence xa', ba' $\in$ E(G) $\Rightarrow$ x'a', b'a' $\notin$ E(G) and a is a unique vertex of both x' and b' of that color. Also b is a unique vertex of a' (else if a'b" $\in$ E(G) where b", b have same color, then <a', b", b, x, b'> = Chair). But then color b by color of y; y by color of a; a', u by color of b, a contradiction. Hence the **Claim** holds.

Finally, let x be a unique vertex of a of that color $\Rightarrow$ a is not a unique vertex of x of that color (else b (y) has another vertex of color of a (x) {else color b, x by color of a; a by color of x; u by color of b} $\Rightarrow$ ba' $\in$ E(G) {else if ba" $\in$ E(G) where a" ($\neq$ a) and a have same color, then <b, a, a", y, a'> = Chair}. Similarly, yx' $\in$ E(G) $\Rightarrow$ x'a' $\notin$ E(G) and <b, a', x', u, x> = Chair). Thus x has multiple vertices of color of a and a unique vertex of color y $\Rightarrow$ xa' $\in$ E(G) (else xa" $\in$ E(G) where a", a have same color and <x, a, a", y, a'> = Chair). Similarly, yx', by', ab' $\in$ E(G) $\Rightarrow$ a'b' $\in$ E(G) (else <x, a', b', u, b> = Chair). Similarly, b'y', y'x', x'a' $\in$ E(G). Now <a', b', y', x'> $\neq$ $K_4$ (else <x', a', b', y', a> = HVN). Hence w.l.g. let a'y' $\notin$ E(G) $\Rightarrow$ $\exists$ y" of color of y s.t. a'y" $\in$ E(G) $\Rightarrow$ x'y" $\notin$ E(G) (else <x', y", y', y, u> = Chair) and <a', x', y", x, u> = Chair, a contradiction.

**Case 2.2**: yb, xa both $\notin$ E(G).

$\forall$ b $\in$ $R_1$, let $T_b$ = {v' $\in$ V(G)-N(u)/ $\exists$ v $\in$ R-$R_1$ s.t. v, v' have same color and bv' $\in$ E(G)}. Let x', y' $\in$ $T_b$ and a', b' $\in$ $T_y$.

**Claim**: <$T_b$> is complete.

Let if possible $\exists$ x', y' $\in$ $T_b$ s.t. x'y' $\notin$ E(G) and let x, y be their corresponding vertices in R-$R_1$. Now x'y (y'x) $\in$ E(G) (else <b, x, y', u, y> = Chair) and clearly, y (x) is a unique vertex of x' (y') of that color. W.l.g. let ax' $\in$ E(G) (else <b, a, x', y', x> = Chair). Now a'b' $\in$ E(G) (else color b by color of x; x' by color of y; y by color of a; a', u by color of b). As G is HVN-free, either xa' or xb' $\notin$ E(G).

**1.** xa' $\notin$ E(G).

$\Rightarrow$ a'x' $\in$ E(G) (else <y, x, a', x', a> = Chair) and $\exists$ a" $\in$ V(G)-N(u) of color of a s.t. xa" $\in$ E(G). Now x'a" $\notin$ E(G) (else <x", a', a", a, u> = Chair) $\Rightarrow$ y'a" $\in$ E(G) (else <x, y', a", y, x'> = Chair) $\Rightarrow$ ba" $\in$ E(G) (else <b, u, x', y', a"> = Chair). Now if ba' $\in$ E(G), then <b, a, a", a', y> = Chair) and if ba' $\notin$ E(G), then <b, u, a", x', a'> = Chair), a contradiction.

**2.** xa' $\in$ E(G) and xb' $\notin$ E(G).

Then x'b' $\in$ E(G) (else <y, x, b', x', b> = Chair) $\Rightarrow$ x'a' $\notin$ E(G) (else <x', b', a', y, x> = HVN). Now x'a $\in$ E(G) (else <y, a', x', u, a> = Chair) $\Rightarrow$ $\exists$ a" $\in$ V(G)-N(u) where a", a have same color s.t. x'a" $\in$ E(G) and <x', a, a", y, a'> = Chair, a contradiction. Hence the **Claim** holds.

As G is HVN-free, w.l.g. let ax' $\notin$ E(G) and P be an a-x bicolor path. Now $\forall$ a $\in$ $R_1$, b has another vertex say a' of color of a (else as G is Chair-free, b is not adjacent to any vertex on P. Alter colors along P, color b by color of a, u by color of b). Let M = {a' $\in$ V(G)-N(u)/ $\exists$ a $\in$ $R_1$ with same color a a and ba' $\in$ E(G)}. Then |M| = |$R_1$|-1. Next we show that <$T_b \cup$ M> is complete. Let if possible $\exists$ x' $\in$ $T_b$, a' $\in$ M s.t. x'a' $\notin$ E(G) $\Rightarrow$ a'x $\in$ E(G) (else <b, a', x', u, x> = Chair) and x is a unique vertex of a' of that color (else if a'x" $\in$ E(G) where x", x have same color, then <a', x, x", b, x'> = Chair). Clearly, a' is a unique vertex of x of that color. But then color x by color of a; a', u by color of x, a contradiction. Hence x'a' $\in$ E(G) $\forall$ x' $\in$ $T_b$ and a' $\in$ M. Finally, we prove that <M> is complete. Let c $\in$ $R_1$ and c', a' $\in$ M. As before cx', cy' $\in$ E(G) $\Rightarrow$ a'c' $\in$ E(G) (else as G is HVN-free, ac', ca' $\notin$ E(G). Now either xa' or xc' $\in$ E(G) {else <b, a', c', u, x> = Chair}. Also xa' $\in$ E(G) iff xc' $\in$ E(G) {else <b, c, c', a', x> = Chair} $\Rightarrow$ <x, a', c', u, a> = Chair). Thus <M> is complete and <b $\cup$ $T_b$ $\cup$ M> =

$K_{|R|} \Rightarrow |R| \leq \omega$. Thus in only possible configuration of $<R>$, every vertex of $R-R_1$ is non-adjacent to every vertex of $R_1$. Also as $<R-R_1>$ is complete, $|R_2| \leq |R_1| \Rightarrow |R_2| \leq \omega/2$.  **II**

**A**: $\exists$ a, b $\in R_1$ with a common $\alpha_i$-vertex say $\alpha_{i1}$ in $<N(u)>$ for some i $\leq \chi$-|R|-1.
Then $c\alpha_{i1} \in E(G)$ $\forall$ c $\in R_1$ (else let $\exists$ c $\in R_1$, s.t. $c\alpha_{i1} \notin E(G)$. By **Lemma**, let $c\alpha_{i2} \in E(G) \Rightarrow$ either $a\alpha_{i2}$ or $b\alpha_{i2} \in E(G)$ {else $<a, b, c, u, \alpha_{i2}>$ = HVN}. W.l.g. let $b\alpha_{i2} \in E(G) \Rightarrow <c, \alpha_{i2}, b, u, \alpha_{i1}>$ = HVN) $\Rightarrow$ a, b have common $\alpha_j$-vertices in N(u) $\forall$ j $\leq \chi$-|R|-1 (else let $\exists$ j s.t. a, b have no common $\alpha_j$-vertex in N(u). Let $a\alpha_{j1}$, $b\alpha_{j2} \in E(G)$. As G is HVN-free, clearly $\alpha_{i1}\alpha_{j1}$, $\alpha_{i1}\alpha_{j2} \in E(G) \Rightarrow <\alpha_{j2}, b, \alpha_{i1}, u, \alpha_{j1}>$ = HVN). Thus $\exists$ T = $\{\alpha_{i1},..., \alpha_{\chi-|R|-1 1}\} \subseteq$ N(u)-R s.t. $<u\cup T\cup R_1>$ is complete (else let $\exists$ $\alpha_{i1}$, $\alpha_{j1}$ common to a and b s.t. $\alpha_{i1}\alpha_{j1} \notin E(G)$. W.l.g. let i < j. Then by **Lemma**, $\alpha_{i1}$ has a $\alpha_j$-vertex say $\alpha_{j2}$ in N(u). As G is HVN-free, w.l.g. let $a\alpha_{j2} \in E(G) \Rightarrow <\alpha_{i1}, \alpha_{j2}, u, a, \alpha_{j1}>$ = HVN).

Thus $1+\chi$-|R|-1+|R_1| $\leq \omega \Rightarrow \chi \leq \omega$+|R|-|R_1| $\leq \frac{3\omega}{2}$ by **II**, a contradiction.

**B**: No a, b $\in R_1$ have a common $\alpha_j$-vertex in $<N(u)>$ for any j.
Let $a\alpha_{i1}$, $b\alpha_{i2} \in E(G)$ $\forall$ i. As G is HVN-free, clearly, $|R_1|$ = 2 and a (b) has no other $\alpha_i$-vertex in $<N(u)>$ for any i. Also $\chi$-|R|-1 > 0 (else $\chi$ = |R|+1 $\leq \omega$+1).
**Case B.1**: $\chi$-|R|-1 > 1.
As G is HVN-free $\alpha_{i1}\alpha_{j1}$ and $\alpha_{i2}\alpha_{j2} \notin E(G)$ and as G is -free, $\alpha_{i2}\alpha_{j1}$ and $\alpha_{i1}\alpha_{j2} \in E(G)$ (i $\neq$ j). Now a (b) has at most one unique $\alpha_j$-vertex (else if a has a unique $\alpha_i$-vertex and $\alpha_j$-vertex, then clearly $\alpha_{i1}$ is a unique $\alpha_i$–vertex of $\alpha_{j2}$ {else if $\alpha_{j2}\alpha_{i3} \in E(G)$, then $<b, a, \alpha_{i2}, \alpha_{j2}, \alpha_{i3}>$ = Chair}. Then color a, $\alpha_{j2}$ by $\alpha_i$; $\alpha_{i1}$ by $\alpha_j$; u by color of a).
**B.1.1**: a has no $\alpha_i$-vertex unique for any i.
Then $\exists$ $\alpha_1$', $\alpha_2$' $\in$ V(G)-N(u) s.t. $a\alpha_1$', $a\alpha_2$' $\in E(G) \Rightarrow b\alpha_1$', $b\alpha_2$' $\in E(G)$ (else $<a, \alpha_1$', $\alpha_{11}, b, \alpha_{12}>$ = Chair) $\Rightarrow \alpha_1$'$\alpha_2$' $\notin E(G)$ (else $<\alpha_1$', $\alpha_2$', a, b, u>$ = HVN). Also $\alpha_1$'$\alpha_{21} \in E(G)$ (else $<a, \alpha_1$', $\alpha_{11}, \alpha_{21}, \alpha_{12}>$ = Chair). Similarly, $\alpha_1$'$\alpha_{22}$, $\alpha_2$'$\alpha_{11}$, $\alpha_2$'$\alpha_{12} \in E(G)$. Again $N_j$ = 2 $\forall$ j (else if $u\alpha_{j3} \in E(G)$, then $<u, \alpha_{j2}, \alpha_{j3}, a, \alpha_j$'>$ = Chair) and $\alpha_j$' is the only $\alpha_j$-vertex of $\alpha_{ik}$ (k = 1, 2) in V(G)-N(u) (else if $\alpha_{ik}\alpha_j$'' is such a vertex, then $<\alpha_{ik}, \alpha_j$'', $\alpha_j$', u, $\alpha_{jk}>$ = Chair). Now $\alpha_{11}$ or $\alpha_{21}$ has a vertex b' of color of b (else color $\alpha_{11}$, $\alpha_{21}$ by color of b; $\alpha_{12}$ by $\alpha_2$; $\alpha_2$', u by $\alpha_1$). Again b'$\alpha_{11} \in E(G)$ iff b'$\alpha_{21} \in E(G)$ (else ab' $\in$ E(G) {else $<a, b, \alpha_{21}, \alpha_{11}, b$'>$ = Chair}, b'$\alpha_2$' $\notin E(G)$ {else $<\alpha_2$', b', a, $\alpha_{11}$, u>$ = HVN} and $<\alpha_{11}, \alpha_2$', b', u, $\alpha_{21}>$ = Chair). Hence b'$\alpha_{11}$, b'$\alpha_{21}$, a'$\alpha_{12}$, a'$\alpha_{22}$ $\in$ E(G).
**B.1.1.1**: Let a'$\alpha_1$' $\in E(G)$.
$\Rightarrow$ ba' $\notin E(G)$ (else $<\alpha_1$', a', $\alpha_{22}$, b, u>$ = HVN) and a is a unique vertex of b of that color. Also a'$\alpha_2$' $\in E(G)$ (else $<b, \alpha_2$', u, $\alpha_1$', a' >$ = Chair)

**Claim**: b'$\alpha_1$' $\notin E(G)$
Let if possible b'$\alpha_1$' $\in E(G)$. As before b'$\alpha_2$' $\in E(G)$. Also ab' $\notin E(G)$ and b is a unique vertex of a of that color $\Rightarrow$ a'b' $\in E(G)$ (else $<\alpha_1$', b', a', a, u>$ = Chair) $\Rightarrow$ b'$\alpha_{12} \notin E(G)$ (else $<\alpha_2$', $\alpha_{12}$, a', b', $\alpha_1$'>$ = HVN). Similarly, b'$\alpha_{22} \notin E(G)$. Again b is a unique vertex of $\alpha_{12}$ of that color (else if $\alpha_{12}$b'' $\in$ E(G) where b'', b have same color, then $<\alpha_{12}, b, b$'', $\alpha_{21}$, b'>$ = Chair). Thus b is a unique vertex of a, $\alpha_{12}$, $\alpha_{22}$ of that color. Color a, $\alpha_{12}$, $\alpha_{22}$ by color of b; b by color of a; $\alpha_{21}$ by $\alpha_1$; $\alpha_1$', u by $\alpha_2$, a contradiction. Hence the **Claim** holds.

Then clearly, b'$\alpha_2$' $\notin E(G)$ and b is a unique vertex of $\alpha_1$'; $\alpha_2$'of that color $\Rightarrow \alpha_{11}$ or $\alpha_{21}$ has another vertex of color of a (else color b, $\alpha_{11}$, $\alpha_{21}$ by color of a; a by $\alpha_1$; $\alpha_1$', u by color of b). W.l.g. let $\alpha_{11}$ have another vertex of color of a $\Rightarrow$ a'$\alpha_{11} \in E(G)$ (else if a''$\alpha_{11} \in E(G)$ where a'', a have same color, then $<\alpha_{11}, a, a$'', $\alpha_{22}$, a'>$ = Chair). Also b'$\alpha_{22} \in E(G)$ (else $<\alpha_{11}, b', \alpha_2$', $\alpha_{22}, \alpha_1$'>$ = Chair) $\Rightarrow$ a'b' $\notin$ E(G) (else $<a', b', \alpha_{11}, \alpha_{22}, u>$ = HVN). But then $<\alpha_{11}, u, b', a', \alpha_1$ '>$ = Chair, a contradiction.
**B.1.1.2**: a'$\alpha_1$', a'$\alpha_2$', b'$\alpha_1$', b'$\alpha_2$' $\notin E(G)$.
$\Rightarrow$ a, b are unique vertices of $\alpha_1$', $\alpha_2$' of those colors (else if a''$\alpha_1$' $\in$ E(G) where a'', a have same color, then $<\alpha_1$', a, a'', $\alpha_{22}$, a'>$ = Chair). Now a' is adjacent to $\alpha_{11}$ or $\alpha_{21}$ (else $<\alpha_{12}, \alpha_{21}, a', \alpha_2$', $\alpha_{11}>$ = Chair). W.l.g. let a'$\alpha_{11} \in E(G) \Rightarrow$ a'$\alpha_{21} \in E(G)$ (else $<\alpha_{12}, \alpha_{21}, \alpha_2$', a', $\alpha_{22}>$ = Chair). Similarly,

b'$\alpha_{12}$, b'$\alpha_{22}$ ∈ E(G) ⇒ a'b' ∉ E(G) (else <b', a', $\alpha_{11}$, $\alpha_{22}$, u> = HVN) ⇒ ab' ∈ E(G) (else <$\alpha_{12}$, a', b', $\alpha_2$', a> = Chair) and a is a unique vertex of that color of b'. But then color a by $\alpha_1$; $\alpha_{11}$ by color of b; b', $\alpha_1$', $\alpha_2$', u by color of a, a contradiction.

**B.1.2**: a (b) has a unique $\alpha_1$-vertex and multiple $\alpha_i$-vertices ∀ i > 1.

⇒ Let a$\alpha_2$' ∈ E(G). As before b$\alpha_2$', $\alpha_2$'$\alpha_{11}$, $\alpha_2$'$\alpha_{12}$ ∈ E(G). Again $N_j$ = 2 ∀ j > 1 (else if u has a $\alpha_{j3}$, then <u, $\alpha_{j2}$, $\alpha_{j3}$, a, $\alpha_j$'> = Chair). Also ∃ a', b' s.t. b'$\alpha_{11}$, b'$\alpha_{21}$, a'$\alpha_{12}$, a'$\alpha_{22}$ ∈ E(G).

**B.1.2.1**: a'$\alpha_2$' ∈ E(G).

⇒ ba' ∉ E(G), a is a unique vertex of b of that color and a'b' ∈ E(G) ⇒ $\alpha_{11}$a' ∈ E(G) (else if a is a unique vertex of $\alpha_{11}$ of that color, then color b, $\alpha_{11}$ by color of a; a by $\alpha_1$; u by color of b and if $\alpha_{11}$a'' ∈ E(G) where a'', a have same color, then <$\alpha_{22}$, b, a', $\alpha_{11}$, a''> = Chair) and a'b' ∈ E(G) (else <$\alpha_{11}$, a', b', u, b> = Chair) ⇒ b'$\alpha_2$' ∉ E(G) (else b'$\alpha_{22}$ ∈ E(G) {else <$\alpha_2$', b', $\alpha_{11}$, a', $\alpha_{22}$> = HVN} and <a', b', $\alpha_{11}$, $\alpha_{22}$, u> = HVN) ⇒ b'$\alpha_{22}$ ∈ E(G) (else <$\alpha_{11}$, $\alpha_{22}$, $\alpha_2$', b', $\alpha_{21}$> = Chair). But then <a', b', $\alpha_{11}$, $\alpha_{22}$, u> = HVN, a contradiction.

**B.1.2.2**: a'$\alpha_2$', b'$\alpha_2$' ∉ E(G).

⇒ a'$\alpha_{21}$ ∈ E(G) (else clearly, $\alpha_{22}$ is a unique $\alpha_2$- vertex of a'. Then color $\alpha_{12}$, $\alpha_{22}$ by color of a; a' by $\alpha_2$; b by $\alpha_1$; u by color of b). Now ab' ∈ E(G) (else b is a unique vertex of a of that color and a'b' ∈ E(G) {else <$\alpha_{21}$, b', a', a, $\alpha_2$'> = Chair} ⇒ b'$\alpha_{12}$ ∉ E(G) {else <a', b', $\alpha_{21}$, $\alpha_{12}$, u> = HVN} and $\alpha_{11}$ is a unique $\alpha_1$-vertex of b'. But then color a, b' by $\alpha_1$; $\alpha_{11}$ by color of b; u by color of a). Also b'$\alpha_{12}$ ∈ E(G) (else if $\alpha_{12}$b'' ∈ E(G) where b'', b have same color, then <$\alpha_{12}$, b'', b, $\alpha_{21}$, b'> = Chair and if b is a unique vertex of $\alpha_{12}$ of that color, then color b, $\alpha_{21}$ by $\alpha_1$; $\alpha_{12}$ by color of b; a by $\alpha_2$; $\alpha_2$', u by color of a). Now a'b' ∈ E(G) (else ∃ a'' of color of a s.t. b'a'' ∈ E(G) and <b', a'', a, $\alpha_{22}$, a'> = Chair). But then <a', b', $\alpha_{21}$, $\alpha_{12}$, u> = HVN, a contradiction.

**Case B.2**: $\chi$-|R|-1 = 1.

⇒ |R| ≥ 3 and |R-$R_1$| ≤ |$R_1$| ⇒ |R| ≤ 4 ⇒ $\omega$ ≤ 4 (else $\chi$ ≤ $\omega$+1). Also if R-$R_1$ ⊆ {x, y}, then by **II**, ax, ay, bx, by ∉ E(G) and a (b) has unique vertices of colors of x and y each. Let ax' ∈ E(G).

**Claim**: x is adjacent to either $\alpha_{11}$ or $\alpha_{12}$.

Let if possible x$\alpha_{11}$, x$\alpha_{12}$ ∉ E(G) and x$\alpha_{13}$ ∈ E(G) for $\alpha_{13}$ ∈ N(u). Let xa', xb' ∈ E(G) ⇒ $\alpha_{12}$a' ∈ E(G) (else <u, a, $\alpha_{12}$, x, a'> = Chair) and as a' is a unique vertex of x of that color, a' is a unique vertex of $\alpha_{12}$ of that color. Now $\alpha_{13}$a' ∈ E(G) (else <u, a, $\alpha_{13}$, $\alpha_{12}$, a'> = Chair). Similarly, b'$\alpha_{11}$, $\alpha_{13}$b' ∈ E(G) ⇒ a'b' ∉ E(G) (else <b', a', x, $\alpha_{13}$, u> = HVN) ⇒ ba' ∈ E(G) (else <x, a', b', u, b> = Chair) and b is a unique vertex of a' of that color (else if a'b'' ∈ E(G) where b, b'' have same color ⇒ <a', b, b'', x, b'> = Chair). Now $\alpha_{12}$ is a unique $\alpha_1$–vertex of b (else if ∃ $\alpha_1$–vertex $\alpha_1$' s.t. b$\alpha_1$' ∈ E(G), then $\alpha_1$' ∈ V(G)-N(u) as a, b have no common $\alpha_1$–vertex in N(u) and <u, $\alpha_{11}$, $\alpha_{13}$, b, $\alpha_1$'> = Chair). But then color b by $\alpha_1$; $\alpha_{12}$ by color of a; a', u by color of b, contradiction. Hence the **Claim** holds.

W.l.g. let x$\alpha_{11}$ ∈ E(G).

**Case B.2.1**: a (b) has a unique $\alpha_1$-vertex $\alpha_{11}$ ($\alpha_{12}$).

⇒ $\alpha_{11}$ ($\alpha_{12}$) has a unique vertex of color of b (a) say b' (a'). As u has maximum no. of unique vertices, clearly a (b) has multiple vertices of color of b (a) and ab' (ba') ∈ E(G) (else if ∃ b'' of color of b s.t. ab'' ∈ E(G), then <a, b'', b, $\alpha_{11}$, b'> = Chair).

**Case B.2.1.1**: x$\alpha_{12}$ ∉ E(G).

⇒ xa' ∈ E(G) (else <u, a, x, $\alpha_{12}$, a'> = Chair), x'$\alpha_{12}$ ∈ E(G) and x' is a unique vertex of $\alpha_{12}$ of that color (else if $\alpha_{12}$x'' ∈ E(G) where x'', x have same color, then <$\alpha_{12}$, x', x'', u, x> = Chair) ⇒ x'a' ∈ E(G) (else if x'a'' ∈ E(G) where a'', a have same color, then <$\alpha_{12}$, a', u, x', a''> = Chair) ⇒ bx' ∉ E(G) (else <x', a', b, $\alpha_{12}$, u> = HVN). Let bx'' ∈ E(G) where x'', x have same color. But then a'x'' ∉ E(G) (else <a', x', x'', x, u> = Chair) ⇒ <b, u, x'', a', x'> = Chair, a contradiction.

**Case B.2.1.2**: x$\alpha_{12}$ ∈ E(G).

Now bx' ∈ E(G) (else x'$\alpha_{11}$ ∈ E(G) {else <a, x', b, $\alpha_{11}$, x> = Chair} ⇒ b'x' ∉ E(G) ⇒ ∃ b'' of color of b s.t. x'b'' ∈ E(G) {else x' has no vertex of color of b. Hence color x' by color of b; a by color of x; u by color of a} and <$\alpha_{11}$, b', u, x', b''> = Chair) ⇒ x'a' ∈ E(G) (else xa' ∈ E(G) {else <b, x', a', u,

x> = Chair} $\Rightarrow \exists$ x" of color of x s.t. a'x" $\in$ E(G) and <a', x", x, b, x'> = Chair). Similarly, x'b' $\in$ E(G) $\Rightarrow$ x'$\alpha_{1i}$ $\notin$ E(G) and x is a unique vertex of $\alpha_{1i}$ for i = 1, 2. Also $\alpha_{11}, \alpha_{12}$ are the only $\alpha_1$–vertices of x (else if z is another $\alpha_1$–vertex of x, then <x, z, $\alpha_{11}, \alpha_{12}$, b> = Chair). But then color $\alpha_{1i}$ by color of x; x by $\alpha_1$; a by $\alpha_1$; u by color of a, a contradiction.

**Case B.2.2**: a, b have multiple $\alpha_1$–vertices.

As a, b have no common $\alpha_1$-vertex in <N(u)>, $\exists$ $\alpha_1$-vertex $\alpha_1$' $\in$ V(G)-N(u) s.t. a$\alpha_1$' $\in$ E(G) $\Rightarrow$ b$\alpha_1$' $\in$ E(G).

**Case B.2.2.1**: x$\alpha_{12}$ $\notin$ E(G).

Let xa', xb' $\in$ E(G). As before $\alpha_{12}$x', $\alpha_{12}$a' $\in$ E(G) and x', a' are unique vertices of $\alpha_{12}$ of those colors $\Rightarrow$ x'a' $\in$ E(G). Also x$\alpha_1$' $\in$ E(G) (else <u, x, $\alpha_{12}$, a, $\alpha_1$'> = Chair). Now bx' $\in$ E(G) (else if bx" $\in$ E(G) where x", x have same color, then $\alpha_{11}$x" $\in$ E(G) {else <b, $\alpha_{12}$, x", a, $\alpha_{11}$> = Chair}. But then if $\alpha_{11}$x' $\notin$ E(G) then <$\alpha_{11}$, x, x", a, x'> = Chair and if $\alpha_{11}$x' $\in$ E(G) then <$\alpha_{11}$, x, x', x", b> = Chair) $\Rightarrow$ x'$\alpha_1$' $\notin$ E(G) and ba' $\notin$ E(G) $\Rightarrow$ a is a unique vertex of b of that color and a'b' $\in$ E(G) (else <x, a', b', u, b> = Chair). Now x is a unique vertex of $\alpha_1$' of that color and $\alpha_1$'a' $\in$ E(G) (else <b, u, $\alpha_1$', x', a'> = Chair). Now b'x' $\in$ E(G) iff b'$\alpha_1$' $\in$ E(G) (else <b, u, x', $\alpha_1$', b'> = Chair) $\Rightarrow$ b'x', b'$\alpha_1$' $\notin$ E(G) (else <x, $\alpha_1$', a', b', x'> = HVN) and $\exists$ x" of color of x s.t. b'x" $\in$ E(G) $\Rightarrow$ a'x" $\notin$ E(G) (else <a', x", x, u> = Chair) and <a', x', $\alpha_1$', b', x"> = Chair, a contradiction.

**Case B.2.2.2**: x$\alpha_{12}$ $\in$ E(G).

Now $N_1$ = 2 (else if $\exists$ $\alpha_{13}$ $\in$ N(u), then <u, $\alpha_{12}$, $\alpha_{13}$, a, $\alpha_1$'> = Chair) $\Rightarrow$ Either $\alpha_{11}$ has a vertex of color of b or $\alpha_{12}$ has a vertex of color of a (else color $\alpha_{11}$ by color of b; $\alpha_{12}$ by color of a; u by $\alpha_1$). W.l.g. let $\alpha_{12}$a' $\in$ E(G) where a', a have same color. Clearly, a' is a unique vertex of $\alpha_{12}$ of that color. Now either x'b $\in$ E(G) or x'$\alpha_{11}$ $\in$ E(G) (else <a, x', b, $\alpha_{11}$, x> = Chair). First let bx' $\notin$ E(G) $\Rightarrow$ x'$\alpha_{11}$ $\in$ E(G). Let bx" $\in$ E(G) $\Rightarrow$ $\alpha_{12}$x" $\in$ E(G) (else <b, a, x", $\alpha_{12}$, x> = Chair) and x"a' $\in$ E(G) (else <$\alpha_{12}$, x", a', u, a> = Chair) $\Rightarrow$ ba' $\notin$ E(G) and a is a unique vertex of b of that color. Now xa' $\in$ E(G) (else <$\alpha_{12}$, x, a', b, a> = Chair) $\Rightarrow$ x" has a" ($\neq$ a) of color of a on x-a bicolor path and <x", a', a", b, a> = Chair, a contradiction. Hence bx' $\in$ E(G) $\Rightarrow$ x'$\alpha_1$' $\notin$ E(G) and x'$\alpha_{11}$ $\in$ E(G) (else <b, x', $\alpha_1$', u, $\alpha_{11}$> = Chair). Similarly, x'$\alpha_{12}$ $\in$ E(G) $\Rightarrow$ x$\alpha_1$' $\in$ E(G) (else <b, x', $\alpha_1$', u, x> = Chair) $\Rightarrow$ x'a' $\in$ E(G) (else x'a" $\in$ E(G) where a", a have same color and <$\alpha_{12}$, u, a', x', a"> = Chair) $\Rightarrow$ ba' $\notin$ E(G) and a is a unique vertex of b of that color. Also $\alpha_1$'a' $\in$ E(G) (else <a, u, $\alpha_1$', x', a'> = Chair) and xa' $\in$ E(G) (else <$\alpha_{12}$, x, a', b, a> = Chair). Let xb' $\in$ E(G) $\Rightarrow$ a'b' $\in$ E(G) (else <x, a', b', u, b> = Chair) $\Rightarrow$ $\alpha_{12}$b' $\notin$ E(G) (else <a', b', x, $\alpha_{12}$, u> = HVN) and b'$\alpha_1$' $\notin$ E(G) (else <$\alpha_1$', b', x, a', $\alpha_{12}$> = HVN) $\Rightarrow$ b'x' $\notin$ E(G) (else <x, u, $\alpha_1$', b', x'> = Chair) $\Rightarrow$ $\exists$ x" of color of x s.t. b'x" $\in$ E(G) on x-b bicolor path. Now $\alpha_1$'x" $\notin$ E(G) (else <$\alpha_1$', x", x, b, x'> = Chair) $\Rightarrow$ <x, $\alpha_1$', u, b', x"> = Chair), a contradiction.

Thus 1+$\chi$-|R|-1+|$R_1$| $\leq$ $\omega$ $\Rightarrow$ $\chi$ $\leq$ $\omega$+|R|-|$R_1$| $\leq$ $\frac{3\omega}{2}$ by **II**, a contradiction.

This proves **Theorem 2**.

**Corollary 4**: If G is {HVN, ($P_3 \cup K_1$)}-free, then $\chi(G) \leq \omega+1$.

Proof: From the proof of **Theorem 2**, it follows that as ($P_3 \cup K_1$) is an induced subgraph of Chair, the only possible configuration of <R> is where <R-$R_1$> is complete and every vertex of R-$R_1$ is non-adjacent to every vertex of $R_1$. Further $\chi \leq \omega+|R-R_1|$.

**Claim**: |R-$R_1$| $\leq$ 1.

Let if possible |R-$R_1$| > 1 and x, y $\in$ R-$R_1$. Let {a, y', a', y} be a-y bicolor path. Then a'x $\in$ E(G) (else <x, y, a', a> = $P_3 \cup K_1$). Let {a, x', a', x} be a-x bicolor path. Then {b, y', b', y}, {b, x', b', x} are b-y, b-x bicolor paths $\Rightarrow$ Now x'y' $\in$ E(G) (else <x', b', y', u> = $P_3 \cup K_1$) $\Rightarrow$ <x', y', a, b, u> = HVN, a contradiction. Hence the **Claim holds**.

$\Rightarrow$ $\chi \leq \omega+1$. This proves **Corollary 4**.

**Corollary 5**: If G is {HVN, $(K_2 U 2K_1)$}-free, then $\chi(G) \leq \omega+1$.
Proof: From the proof of **Theorem 2**, it follows that as $(K_2 U 2K_1)$ is an induced subgraph of Chair, the only possible configuration of $<R>$ is where $<R-R_1>$ is complete and every vertex of $R-R_1$ is non-adjacent to a vertex of $R_1$. Further $\chi \leq \omega+|R-R_1|$.

**Claim**: $|R-R_1| \leq 1$.
Let if possible $|R-R_1| > 1$ and $x, y \in R-R_1$. Then as G is $(K_2 U 2K_1)$-free, if bx', by' $\in$ E(G), then ax', ay' $\in$ E(G). Similarly, xa', xb' $\in$ E(G) $\Rightarrow$ As G is HVN-free x'y', a'b' $\notin$ E(G) $\Rightarrow$ x'y, y'x, a'b, b'a $\in$ E(G). But then color a by color of y, y' by color of x; x by color of b; b' and u by color of a, a contradiction. Hence the **Claim** holds.

$\Rightarrow \chi \leq \omega+1$. This proves **Corollary 5**.

**Corollary 6**: If G is {Chair, $K_4$}-free, then $\chi(G) \leq \omega+1$.
Proof: From the proof of **Theorem 2**, it follows that as $(K_2 U 2K_1)$ is an induced subgraph of Chair, the only possible configuration of $<R>$ is where $<R-R_1>$ is complete and every vertex of $R-R_1$ is non-adjacent to a vertex of $R_1$. Further $\chi \leq \omega+|R-R_1|$.

**Claim**: $|R-R_1| \leq 1$.
Let if possible $|R-R_1| > 1$ and $x, y \in R-R_1$. Let bx', by', xa', xb' $\in$ E(G) $\Rightarrow$ Either x'y' or a'b' $\in$ E(G) (else as G is Chair-free, x'y {y'x, a'b, ba'} $\in$ E(G) and y {x, b, a} is a unique vertex of x' {y', a', b'} of that color. Then color b by color of y, y' by color of x; x by color of a; a' and u by color of b). W.l.g. let x'y' $\in$ E(G). As G is $K_4$-free, a is non-adjacent to either x' or y'. W.l.g. let ax' $\notin$ E(G). Then b has another vertex say a' of color of a (else let M be a bicolor component containing a, x. As G is Chair-free, M is a path and b is not adjacent to any vertex of M. Alter colors along M, color b by color of a; u by color of b). Now a'x' $\in$ E(G) (else as G is Chair-free, a'x $\in$ E(G) and x (a') is a unique vertex of a' (x) of that color. Color x by color of a; a', u by color of x). Similarly, a'y' $\in$ E(G) and $<b, x', y', a'> = K_4$, a contradiction. Hence the **Claim** holds.

$\Rightarrow \chi \leq \omega+1$. This proves **Corollary 6**.

**Theorem 3**: If G is {$K_5$-e, $(P_3 U K_1)$}-free, then $\chi(G) \leq \omega+1$.
Proof: Let if possible G be a smallest such graph with $\chi(G) > \omega + 1$. Among all $(\chi-1)$-colorings of G-v for every v $\in$ V(G), let $\exists$ u $\in$ V(G) and $(\chi-1)$-coloring C of G-u s.t. u has maximum no. |R| vertices with unique colors and if u has repeat colors $\alpha_1, ..., \alpha_{\chi-|R|-1}$, then u has minimum no. $N_1$ of $\alpha_1$-vertices. Further with |R| unique colors and $N_1$ no. of $\alpha_1$-vertices, let u have minimum no. $N_2$ of $\alpha_2$-vertices and so on. Let R = {x $\in$ N(u)/ x receives a unique color from C in $<N(u)>$}. Let $<R_1>$ be a maximum clique in $<R>$. Let S = {x' $\in$ V(G)-N(u)/ $\exists$ x $\in$ R_1 with same color s.t. xa $\notin$ R_1 for some a $\in$ R_1}. Let T = {x $\in$ R/ xa $\in$ E(G) $\forall$ a $\in$ R} and T' = {x $\in$ T/ $\exists$ s' $\in$ S s.t. xs' $\notin$ E(G)}.

**Claim**: $|R| \leq \omega$.
First we prove that $<S>$ is complete. Let if possible $\exists$ v', w' $\in$ S, s.t. v'w' $\notin$ E(G). Let v, w be their corresponding vertices in R. Now $\exists$ z $\in$ R s.t. vz $\notin$ E(G) $\Rightarrow$ zw' $\in$ E(G) (else $<u, z, v', w'> = P_3 U K_1$). Similarly, zv' $\in$ E(G) $\Rightarrow <v', z, w', v> = P_3 U K_1$, a contradiction. Hence $<S>$ is complete. Next let x $\in$ T' and s' $\in$ S s.t. xs' $\notin$ E(G). Now by construction of S, $\exists$ s, y $\in$ R_1 s.t. s, s' have same color and sy $\notin$ E(G) $\Rightarrow$ s' is a unique vertex of y of that color and hence s' has a vertex say x' $\in$ V(G)-N(u) of color of x. Let M = {x' $\in$ V(G)-N(u)/ x $\in$ T'}. Clearly, $<S \cup T-T' \cup M> = K_{|R|}$. Thus the **Claim** holds.

$\Rightarrow \chi-|R|-1 > 0$ (else $\chi \leq |R|+1 \leq \omega+1$).     **I**
Let if possible $|R-R_1| > 1$ and $x, y \in R-R_1$.     **II**

**Case 1**: $\exists$ x, y $\in$ |R-R_1| s.t. xy $\notin$ E(G).
**Case 1.1**: $\exists$ a $\in$ R_1, s.t. ax, ay $\notin$ E(G).

Let xa' $\in$ E(G) $\Rightarrow$ ya' $\in$ E(G) (else <a, u, y, a'> = $P_3 \cup K_1$). But then <x, a', y, a> = $P_3 \cup K_1$, a contradiction.

**Case 1.2**: **Case 1.1** does not exist.

Let ay $\notin$ E(G) $\Rightarrow$ ax $\in$ E(G). Let b $\in$ $R_1$ be s.t. xb $\notin$ E(G) $\Rightarrow$ yb $\in$ E(G). Let {x, y', x', y} be a bicolor x-y path $\Rightarrow$ ay' $\in$ E(G) (else <a, x, y', y> = $P_3 \cup K_1$). Similarly, bx' $\in$ E(G). Let xb', ya' $\in$ E(G). As G is ($P_3 \cup K_1$)-free, clearly, a'y', x'b' $\in$ E(G) and <a', b', x', y'> = $K_4$ $\Rightarrow$ yb' $\notin$ E(G) (else <y', b', a', x', y> = $K_5$-e) and b is a unique vertex of y of that color. Similarly, xa' $\notin$ E(G). Again ab' $\in$ E(G) (else <a, y', b', y> = $P_3 \cup K_1$). Similarly, ba' $\in$ E(G) $\Rightarrow$ As G is ($K_5$-e)-free ax', by' $\notin$ E(G). By **Lemma**, let a$\alpha_{11}$ $\in$ E(G) where $\alpha_{11}$ $\in$ N(u)-R. Then b$\alpha_{11}$ $\notin$ E(G) (else as G is ($K_5$-e)-free, x$\alpha_{11}$, y$\alpha_{11}$ $\notin$ E(G) $\Rightarrow$ x'$\alpha_{11}$ $\in$ E(G) (else <x, u, $\alpha_{11}$, x'> = $P_3 \cup K_1$), a'$\alpha_{11}$ $\in$ E(G) (else <x, u, $\alpha_{11}$, a'> = $P_3 \cup K_1$) and <y, a', x', b, $\alpha_{11}$> = ($K_5$-e)) $\Rightarrow$ b'$\alpha_{11}$ $\in$ E(G) (else <b, u, $\alpha_{11}$, b'> = $P_3 \cup K_1$). Also y'$\alpha_{11}$ $\in$ E(G) (else <b, u, $\alpha_{11}$, y'> = $P_3 \cup K_1$) $\Rightarrow$ x$\alpha_{11}$ $\in$ E(G) (else <x, a, b', y', $\alpha_{11}$> = ($K_5$-e)). But then <u, x, a, $\alpha_{11}$, b'> = ($K_5$-e), a contradiction.

**Case 2**: <R-$R_1$> is complete.

W.l.g. bx, ay $\notin$ E(G). Then either yb, xa both $\notin$ E(G) or both $\in$ E(G) (else this case can be proved as **Case 1.2**).

**Case 2.1**: yb, xa both $\in$ E(G).

Now |$R_1$| = 2 (else if $\exists$ c $\in$ $R_1$, then as G is $K_5$-e, cx, cy $\notin$ E(G) and with <c, a, x, y> this becomes **Case 1.2**) $\Rightarrow$ R-$R_1$ = 2 and |R| = 4. Let {a, y', a', y} and {b, x', b', x} be a-y and b-x bicolor paths. Now a' is adjacent to b or x (else <b, u, x, a'> = $P_3 \cup K_1$). W.l.g. let a'b $\in$ E(G) $\Rightarrow$ a'b' $\in$ E(G) (else <a', b, u, b'> = $P_3 \cup K_1$). As G is ($P_3 \cup K_1$)-free, <a', b', x', y'> = $K_4$ $\Rightarrow$ by' $\notin$ E(G) (else <b', x', y', a', b> = $K_5$-e) and y'x $\in$ E(G) (else <b, u, x, y'> = $P_3 \cup K_1$) $\Rightarrow$ xa' $\notin$ E(G). Also as $\chi$-|R|-1 > 0 by **Lemma**, let b$\alpha_{11}$ $\in$ E(G) where $\alpha_{11}$ $\in$ N(u)-R. Now y$\alpha_{11}$ $\notin$ E(G) (else clearly, a$\alpha_{11}$ $\notin$ E(G) $\Rightarrow$ a'$\alpha_{11}$ $\in$ E(G) (else <a, u, $\alpha_{11}$, a'> $P_3 \cup K_1$) and <u, b, $\alpha_{11}$, y, a'> = $K_5$-e) $\Rightarrow$ y'$\alpha_{11}$ $\in$ E(G) and a$\alpha_{11}$ $\in$ E(G) (else <a, y', $\alpha_{11}$, y> $P_3 \cup K_1$). Also x$\alpha_{11}$ $\notin$ E(G) (else <b, u, a, $\alpha_{11}$, x> = $K_5$-e) $\Rightarrow$ x'$\alpha_{11}$ $\in$ E(G) (else <$\alpha_{11}$, u, x, x'> = $P_3 \cup K_1$) and a'$\alpha_{11}$ $\in$ E(G) (else <$\alpha_{11}$, u, x, a'> = $P_3 \cup K_1$) $\Rightarrow$ x'a $\notin$ E(G) (else <u, b, a, $\alpha_{11}$, x'> = $K_5$-e) and x'y $\in$ E(G) (else <a, u, y, x'> = $P_3 \cup K_1$). But then <y, b, x', a', $\alpha_{11}$> $K_5$-e, a contradiction.

**Case 2.2**: yb, xa both $\notin$ E(G).

Let {a, y', a', y} and {a, x', a', x} be a-y and a-x bicolor paths. If xb' $\in$ E(G), then as G is $P_3 \cup K_1$-free, {b, y', b', y} and {b, x', b', x} are b-y and b-x bicolor paths. Also as G is $P_3 \cup K_1$-free, <a', b', x', y'> = $K_4$ $\Rightarrow$ As G is ($K_5$-e)-free xy', yx', ab', ba' $\notin$ E(G). As G is ($K_5$-e)-free, |$R_1$| = 2 and |R| = 4. Let x'$\alpha_{11}$ $\in$ E(G) with $\alpha_{11}$ $\in$ N(u)-R (else color x' by $\alpha_1$; a by color of x; u by color of a). Now b$\alpha_{11}$ $\notin$ E(G) (else a$\alpha_{11}$ $\notin$ E(G) {else <u, a, b, $\alpha_{11}$, x'> = $K_5$-e} $\Rightarrow$ a'$\alpha_{11}$ $\in$ E(G) {else <a, u, $\alpha_{11}$, a'> = $P_3 \cup K_1$}, x$\alpha_{11}$ $\in$ E(G) {else <a, x', $\alpha_{11}$, x> = $P_3 \cup K_1$}, y$\alpha_{11}$ $\in$ E(G) and <u, x, y, $\alpha_{11}$, a'> = $K_5$-e, a contradiction. Similarly, a$\alpha_{11}$ $\notin$ E(G) $\Rightarrow$ a'$\alpha_{11}$, b'$\alpha_{11}$ $\in$ E(G) and x$\alpha_{11}$ $\in$ E(G) (else <b, x', $\alpha_{11}$, x> = $P_3 \cup K_1$). But then <x, $\alpha_{11}$, a', b', x'> = $K_5$-e, a contradiction.

Hence our assumption **II** is wrong and |R-$R_1$| $\le$ 1                                                          **III**

**A**: $\exists$ a, b $\in$ $R_1$ with a common $\alpha_i$-vertex say $\alpha_{i1}$ for some i, 1 $\le$ i $\le$ $\chi$-|R|-1.

Clearly as G is ($K_5$-e)-free, $\alpha_{i1}$ is adjacent to all vertices of $R_1$. Now a, b have a common $\alpha_j$-vertex $\forall$ j (else by **Lemma** a$\alpha_{j1}$, b$\alpha_{j2}$ $\in$ E(G). As G is ($K_5$-e)-free, $\alpha_{i1}\alpha_{j1}$, $\alpha_{i1}\alpha_{j2}$ $\notin$ E(G) $\Rightarrow$ <$\alpha_{i1}$, a, $\alpha_{j1}$, $\alpha_{j2}$> = $P_3 \cup K_1$). Hence all vertices of $R_1$ have a common $\alpha_i$-vertex $\forall$ i. Also as G is ($K_5$-e)-free, $\exists$ T = {$\alpha_{i1}$, …, $\alpha_{\chi-|R|-1 1}$} s.t. <$R_1 \cup u \cup T$> = $K_{|R_1|+1+\chi-|R|-1}$ $\Rightarrow$ $\chi \le \omega+1$ by **III**, a contradiction.

**B**: No a, b $\in$ $R_1$ have a common $\alpha_i$-vertex in N(u) for any i $\le$ $\chi$-|R|-1.

Let a$\alpha_{i1}$, b$\alpha_{i2}$ $\in$ E(G). As G is $P_3 \cup K_1$-free, $\alpha_{i1}$ ($\alpha_{i2}$) is the unique $\alpha_i$-vertex of a (b) and G has only these two $\alpha_i$-vertices (else if $\exists$ another $\alpha_i$-vertex $\alpha_{i3}$, then <b, a, $\alpha_{i1}$, $\alpha_{i3}$> = $P_3 \cup K_1$). Now |$R_1$| = 2 (else if $\exists$ c $\in$ $R_1$ then c$\alpha_{i1}$, c$\alpha_{i2}$ $\notin$ E(G) and c has no $\alpha_i$-vertex in N(u)) $\Rightarrow$ |R| $\le$ 3. Now $\chi$-|R|-1 $\le$ 1 (else even if |R| = 3 and x $\in$ R-$R_1$, then w.l.g. let xa $\notin$ E(G). Let u be colored with color of a and a = u'. Thus as a has two vertices say $\alpha_{i1}$, $\alpha_{j1}$ of unique colors clearly, a has |R|+1 vertices of unique colors in

the new coloring) $\Rightarrow \chi-|R|-1 = 1$ and $|R| = 3 = \omega$ (else $\chi \leq \omega+1$). Let $x \in R-R_1$. By **Lemma**, w.l.g. let $x\alpha_{11} \in E(G)$. Then $xa \notin E(G)$ (else $\omega \geq 4$) $\Rightarrow x\alpha_{12} \in E(G)$ (else $<a, \alpha_{11}, x, \alpha_{12}> = P_3 U K_1$) and $xb \notin E(G)$. Let $\{a, x', a', x\}$ be a bicolor a-x path $\Rightarrow bx' \in E(G)$ (else $<x', a, b, x> = P_3 U K_1$). Then $a'\alpha_{12} \in E(G)$ (else $<\alpha_{12}, u, a, a'> = P_3 U K_1$). Let $\{b, x', b', x\}$ be a bicolor b-x path. Then as before $b'\alpha_{11} \in E(G)$. Also $a'b' \in E(G)$ (else $<a', x', b', u> = P_3 U K_1$) $\Rightarrow a'\alpha_{11}, b'\alpha_{12} \notin E(G)$ (else $\omega \geq 4$) $\Rightarrow ab' \in E(G)$ (else $<\alpha_{12}, b, a, b'> = P_3 U K_1$). Similarly, $ba' \in E(G) \Rightarrow x'\alpha_{11}, x'\alpha_{12} \notin E(G)$ (else $\omega \geq 4$) and $<\alpha_{11}, u, \alpha_{12}, x'> = P_3 U K_1$ a contradiction.

This proves **Theorem 3**.

**Theorem 4**: If G is $\{K_5-e, (K_2 U 2K_1)\}$-free, then $\chi(G) \leq \omega+1$.
Proof: Let if possible G be a smallest such graph with $\chi(G) > \omega + 1$. Among all $(\chi-1)$-colorings of G-v for every $v \in V(G)$, let $\exists u \in V(G)$ and $(\chi-1)$-coloring C of G-u s.t. u has maximum no. $|R|$ vertices with unique colors and if u has repeat colors $\alpha_1, ..., \alpha_{\chi-|R|-1}$, then u has minimum no. $N_1$ of $\alpha_1$-vertices. Further with $|R|$ unique colors and $N_1$ no. of $\alpha_1$-vertices, let u have minimum no. $N_2$ of $\alpha_2$-vertices and so on. Let $R = \{x \in N(u) / x$ receives a unique color from C in $<N(u)>\}$. Let $<R_1>$ be a maximum clique in $<R>$. Let if possible $|R-R_1| > 1$ and $x, y \in R-R_1$. **I**

**Case 1**: $\exists x, y \in R-R_1$ s.t. $xy \notin E(G)$.
**Case 1.1**: $\exists a \in R_1$ s.t. $ax, ay \notin E(G)$.
Let $\{a, y', a', y\}, \{a, x', a', x\}$ be a-x, a-y bicolor paths. As G is $(K_2 U 2K_1)$-free, $xy', x'y', yx' \in E(G)$.

**Claim**: $\exists v \in V(G)$ s.t. v has $|R|$ vertices with unique colors, $N_i$ no. of $\alpha_i$-vertices $\forall i$ and if R' is a set of vertices of v with unique colors and $<R_1'>$ is a maximum clique in $<R'>$, then $|R_1'| > 1$.
Let $|R_1| = 1 \Rightarrow <R> = \overline{K}_{|R|}$. As G is $(K_5-e)$-free, $|R-R_1| = 2$ (else if $\exists z \in R-R_1$, then let $z' \in V(G)-N(u)$ be s.t. z, z' have same color and $z'x \in E(G) \Rightarrow <a, x', y', z', a'> = (K_5-e)$). Let a have $N_i'$ no. of $\alpha_i$-vertices. If $N_i' = N_i \forall i$, then color u by color of a and $a = u'$, then in this coloring C', $\{u, x', y'\} \subseteq R' = \{v \in N(a) / v$ receives a unique color in $<N(a)>\}$ and by our assumption $|R'| = |R|$ vertices. Let $R_1'$ be a maximum clique in $<R'> \Rightarrow |R_1'| \geq 2$ and this proves the **Claim**.

Next let $N_i' > N_i$ for some i. As G is $(K_2 U 2K_1)$-free, $\exists$ at most one $\alpha_i$-vertex say $\alpha_i' \in V(G)-N(u)$ and $a\alpha_{ij}, a\alpha_i' \in E(G) \forall \alpha_{ij} \in N(u)$. As G is $(K_2 U 2K_1)$-free, $\alpha_i'$ is adjacent to at least two of $\{x', y', a'\}$ and as G is $(K_5-e)$-free, $\alpha_i'$ is non-adjacent to at least one of $\{x', y', a'\}$. W.l.g. let $\alpha_i'y' \in E(G)$.
- $\alpha_i'x' \in E(G)$.
$\Rightarrow \alpha_i'a' \notin E(G) \Rightarrow \alpha_{i1}a', \alpha_{i2}a' \in E(G)$. Now x', y' have no common $\alpha_i$-vertex in N(u) (else if w is such a vertex, then as G is $(K_2 U 2K_1)$-free, $<w, a, x', y', \alpha_i'> = K_5-e$) $\Rightarrow |N_i| = 2$ (else x, y have a common $\alpha_i$-vertex in N(u)). W.l.g. let $x'\alpha_{i1}, y'\alpha_{i2} \notin E(G)$ and $x'\alpha_{i2}, y'\alpha_{i1} \in E(G) \Rightarrow$ Now $x\alpha_{i1} \in E(G)$ (else color $\alpha_{i1}$ by color of x; y' by $\alpha_i, \alpha_i'$, u by color of a). Similarly, $y\alpha_{i2} \in E(G)$. But then color $\alpha_{i1}$ by color of x; x, y by color of a; $\alpha_{i2}$ by color of y; a', u by $\alpha_i$, a contradiction.
- $\alpha_i'x' \notin E(G)$.
$\Rightarrow \alpha_i'a' \in E(G)$ and as G is $(K_2 U 2K_1)$-free $\alpha_{i1}x', \alpha_{i2}x' \in E(G)$. W.l.g. let $\alpha_{i2}a' \in E(G) \Rightarrow \alpha_{i2}y' \notin E(G)$ (else $<a, \alpha_{i2}, x', y', a'> = K_5-e$) $\Rightarrow y'\alpha_{i1} \in E(G) \Rightarrow \alpha_{i1}a' \notin E(G)$ (else $<a, \alpha_{i1}, x', y', a'> = K_5-e$). Now $\alpha_{i2}y \in E(G)$ (else color $\alpha_{i2}$ by color of y; x by color of a; a' by $\alpha_i$; $\alpha_i'$, u by color of x). But then color $\alpha_{i2}$ by color of y; $\alpha_i'$ by color of x; x, y by color of a; a', u by $\alpha_i$, a contradiction. Hence the **Claim** holds.

Hence w.l.g. let $|R_1| > 1$ and $b \in R_1$. Now b is non-adjacent either to x' or y' or a' (else $<a', x', y', b, a> = K_5-e$). Let $z' \in \{x', y', a'\}$ s.t. $bz' \notin E(G)$ and $z \in \{x, y, a\}$ s.t. z', z have same color. Also let $w \in \{x, y, a\}-z, w' \in \{x', y', a'\}$ s.t. w, w' have same color. Then z' has a unique vertex b' of color of b and $<a', x', y', b'> = K_4$ (else color b' by color of w; w by color of z; z', u by color of b) $\Rightarrow$ b' is non-adjacent to x, y and a and b is a unique vertex of a, x, y of that color. But then color w, z by color of b; b by color of z; u by color of w, a contradiction.

**Case 1.2**: **Case 1.1** does not exist.

Let ay $\notin$ E(G) $\Rightarrow$ ax $\in$ E(G). Let b $\in$ $R_1$ be s.t. xb $\notin$ E(G) $\Rightarrow$ yb $\in$ E(G). As G is ($K_2 \cup 2K_1$)-free, clearly $\exists$ x', y' a', b' s.t. xy', yx', x'y', ay', y'a', a'y, bx', x'b', b'x $\in$ E(G). Now x'a' $\in$ E(G) (else a is a unique vertex of x' of that color $\Rightarrow$ ab' $\in$ E(G) {else color a by color of b; b by color of x; x', u by color of a} $\Rightarrow$ y'b' $\notin$ E(G) {else <x', y', a, b', x> = $K_5$-e}. But then color b by color of x; x' by color of a; a by color of y; y', u by color of b). Similarly, y'b' $\in$ E(G). Again a'b' $\in$ E(G) (else a is a unique vertex of b' of that color $\Rightarrow$ y'b $\in$ E(G). Similarly, x'a $\in$ E(G) and <b, a, y', x', b'> = $K_5$-e). Thus <a', b', x', y'> = $K_4$ $\Rightarrow$ xa', yb' $\notin$ E(G). Also ab' $\in$ E(G) (else b is a unique vertex of a and y of that color $\Rightarrow$ ba', by' $\in$ E(G) and <b', x', y', a', b> = $K_5$-e). Similarly, ba' $\in$ E(G) $\Rightarrow$ ax', by' $\notin$ E(G).

**Claim 1**: $|R| \leq 5$, if $|R| = 4$, then <R> = $P_4$ and if $|R| = 5$, then <R> = $C_5$.

Now $|R_1| = 2$ (else if $\exists$ c $\in$ $R_1$, then as G is ($K_5$-e)-free, cx, cy $\notin$ E(G) and with <u, c, x, y> this can be proved as **Case 1.1**). Also let $|R-R_1| > 2$ and z $\in$ R-$R_1$. Again z is adjacent to either x or y. W.l.g. let zx $\in$ E(G) $\Rightarrow$ As $|R_1| = 2$, za $\notin$ E(G), zy $\in$ E(G), zb $\notin$ E(G) $\Rightarrow$ <R> = $C_5$. This proves **Claim 1**.

Thus $\chi - |R| - 1 > 0$ (else if $|R| = 4$, then as proved above <a', b', x', y'> = $K_4$ and if $|R| = 5$, then again <a', b', x', y', z'> = $K_5$ $\Rightarrow$ $\chi = |R|+1 \leq \omega+1$). By **Lemma**, let $a\alpha_{11} \in$ E(G).
a) $b\alpha_{11} \notin$ E(G) $\Rightarrow$ $\alpha_{11}$ is adjacent to either y' or b' and y' or x $\Rightarrow$ $\alpha_{11}$y' $\in$ E(G) (else $\alpha_{11}$b', $\alpha_{11}$x $\in$ E(G) and <y', x, a, b', $\alpha_{11}$> = $K_5$-e). Again $\alpha_{11}$ is adjacent to either x or b' (else <x, b', $\alpha_{11}$, b> = ($K_2 \cup 2K_1$)) $\Rightarrow$ $\alpha_{11}$ is adjacent to both x and b' (else <x/b', a, y', b'/x, $\alpha_{11}$> = $K_5$-e). But then <u, x, a, $\alpha_{11}$, b'> = $K_5$-e, a contradiction.
b) $b\alpha_{11} \in$ E(G) $\Rightarrow$ As G is $K_5$-e, $x\alpha_{11}$, $y\alpha_{11}$ $\notin$ E(G) $\Rightarrow$ $\alpha_{11}$x' $\in$ E(G) (else <y, x', x, $\alpha_{11}$> = ($K_2 \cup 2K_1$)). Similarly, $\alpha_{11}$y' $\in$ E(G) $\Rightarrow$ b'$\alpha_{11}$ $\notin$ E(G) (else <x, a, y', b', $\alpha_{11}$> = $K_5$-e). Similarly, $\alpha_{11}$a' $\notin$ E(G). But then <x, b', $\alpha_{11}$, y> = ($K_2 \cup 2K_1$), a contradiction.

**Case 2**: <R-$R_1$> is complete.

W.l.g. bx, ay $\notin$ E(G). Then either yb, xa both $\notin$ E(G) or both $\in$ E(G) (else this case can be proved as **Case 1.2**). Also $|R-R_1| \leq |R_1| = 2$ (else if $\exists$ c $\in$ $R_1$, then as G is ($K_5$-e)-free cx,cy $\notin$ E(G) and this case can be proved as **Case 1.2**). Thus $|R| \leq 4$. Let {a, y', a', y} and {b, x', b', x} be a-y and b-x bicolor paths.

**Case 2.1**: yb, xa both $\in$ E(G).

Let {a, y', a', y}, {b, x', b' x} be a-y, b-x bicolor paths.

**Claim**: <a', b', x', y'> = $K_4$.

W.l.g. let x'a' $\notin$ E(G) $\Rightarrow$ ab' $\in$ E(G) {else color a by color of b; b by color of x; x', u by color of a}. Similarly, xy' $\in$ E(G). Now a'b' $\in$ E(G) {else <x', b', a', u> = ($K_2 \cup 2K_1$)}. Similarly, x'y' $\in$ E(G) $\Rightarrow$ y'b' $\notin$ E(G) {else <a, x, y', b', a'> = $K_5$-e}. But then color b by color of x; x' by color of a; a by color of y; y', u by color of b), a contradiction. Hence the **Claim** holds.

Also b is non-adjacent to either a' or y' (else <b', a', x', y', b> = $K_5$-e). W.l.g. let ba' $\notin$ E(G) $\Rightarrow$ xa' $\in$ E(G) (else a is a unique vertex of x, b of that color $\Rightarrow$ ax', ab' $\in$ E(G) and <a, x', y', b', a'> = $K_5$-e) $\Rightarrow$ xy' $\notin$ E(G) (else <x, y', b', a', x'> = $K_5$-e). Now by' $\in$ E(G) (else y is a unique vertex of x, b of that color $\Rightarrow$ yx', yb' $\in$ E(G) and <y, a', x', b', y'> = $K_5$-e). Again a (y) is non-adjacent to either b' or x'. Also either a or y is adjacent to either b' or x' (else <b', x', a, y> = ($K_2 \cup 2K_1$)). W.l.g. let ax' $\in$ E(G) $\Rightarrow$ ab' $\notin$ E(G) $\Rightarrow$ yb' $\in$ E(G) (else color y, a by color of b; b by color of a; u by color of y) $\Rightarrow$ yx' $\notin$ E(G). Now $\chi - |R| - 1 > 0$ (else $\chi = |R|+1 \leq 5 \leq \omega+1$). Hence by **Lemma** let $a\alpha_{11} \in$ E(G) $\Rightarrow$ $b\alpha_{11} \notin$ E(G) (else $x\alpha_{11}$, $y\alpha_{11}$, $x'\alpha_{11}$, $y'\alpha_{11}$ $\notin$ E(G) and <x, y, $\alpha_{11}$, x'> = ($K_2 \cup 2K_1$)). Similarly, $x\alpha_{11} \notin$ E(G) $\Rightarrow$ $x'\alpha_{11} \in$ E(G) (else <b, x', $\alpha_{11}$, x> = ($K_2 \cup 2K_1$)). Similarly, $b'\alpha_{11} \in$ E(G). Now $y'\alpha_{11} \notin$ E(G) (else <b, a, y', x', $\alpha_{11}$> = $K_5$-e) $\Rightarrow$ <y', b, $\alpha_{11}$, x> = ($K_2 \cup 2K_1$), a contradiction.

**Case 2.2**: x, y non-adjacent to both a b.

Let {a, x', a', x}, {b, y', b', y} be a-x, a-y bicolor paths. As G is $(K_2 \cup 2K_1)$-free, clearly {a, y', a', y}, {b, x', b', x} are a-y, b-x bicolor paths. Now x'y' ∈ E(G) or a'b' ∈ E(G) (else color a by color of x; x' by color of y; y by color of b; b', u by color of a). W.l.g. let a'b' ∈ E(G) ⇒ x'y' ∈ E(G) (else x'y ∈ E(G) ⇒ <x, a', b', y', x'> = $K_5$-e) ⇒ ba', ab', xy', yx' ∉ E(G). As χ-|R|-1 > 0, x'$\alpha_{11}$ ∈ E(G) for some $\alpha_{11}$ ∈ N(u) (else color x' by $\alpha_1$; a by color of x; u by color of a) ⇒ $\alpha_{11}$ is non-adjacent to a or b (else <u, a, b, $\alpha_{11}$, x'> = $K_5$-e). W.l.g. let $\alpha_{11}$b ∉ E(G). Then $\alpha_{11}$ is adjacent to a' or b' (else <a', b', b, $\alpha_{11}$> = ($K_2 \cup 2K_1$)) ⇒ $\alpha_{11}$ is non-adjacent to x or y (else <u, $\alpha_{11}$, x, y, a'/b'> = $K_5$-e). W.l.g. let $\alpha_{11}$x ∉ E(G). ⇒ $\alpha_{11}$a' ∈ E(G) (else <x, a', $\alpha_{11}$, b> = ($K_2 \cup 2K_1$)). Similarly, $\alpha_{11}$b' ∈ E(G) ⇒ $\alpha_{11}$y ∉ E(G) (else <$\alpha_{11}$, a', b', y, x> = $K_5$-e). But then <x, y, $\alpha_{11}$, b> = ($K_2 \cup 2K_1$), a contradiction.

Hence **I** is wrong and |R-$R_1$| = 1.                                                                                           **II**

**Case A**: For every i ∈ {1, 2,..., χ-|R|-1}, ∃ some a, b ∈ $R_1$ with a common $\alpha_i$–vertex in <N(u)>.
Clearly, if a, b ∈ $R_1$ have a common $\alpha_i$–vertex say $\alpha_{ij}$ in <N(u)>, then as G is ($K_5$-e)-free, all vertices of $R_1$ are adjacent to $\alpha_{ij}$. Further if a, b have a common $\alpha_i$–vertex $\alpha_{ij}$ and $\alpha_k$–vertex $\alpha_{km}$ in <N(u)>, then as G is ($K_5$-e)-free, $\alpha_{ij}\alpha_k$ ∈ E(G). Thus ∃ T = {$\alpha_{ij}$ ∈ N(u)/ a$\alpha_{ij}$, b$\alpha_{ij}$ ∈ E(G), i, j = 1,..., χ-|R|-1} and <u∪T∪$R_1$> is complete ⇒ 1+χ-|R|-1+|$R_1$| ≤ ω ⇒ χ ≤ ω+1 by **II**, a contradiction.

**Case B:** ∃ i ∈ {1, 2,..., χ-|R|-1}, s.t. no a, b ∈ $R_1$ have a common $\alpha_i$–vertex in <N(u)>.
Let a$\alpha_{i1}$, b$\alpha_{i2}$ ∈ E(G) ⇒ $N_i$ = 2 (else if ∃ $\alpha_{i3}$ ∈ N(u), then as G is ($K_2 \cup 2K_1$)-free, a$\alpha_{i3}$, b$\alpha_{i3}$ ∈ E(G)). Also |$R_1$| = 2 (else if ∃ c ∈ $R_1$, then as G is ($K_5$-e)-free, c$\alpha_{i1}$, c$\alpha_{i2}$ ∉ E(G), contrary to **Lemma**).

**Claim**: No a, b ∈ $R_1$ have a common $\alpha_j$–vertex in <N(u)> for any j ∈ {1, 2,..., χ-|R|-1}-{i}.
Let if possible a, b have a common $\alpha_j$–vertex say $\alpha_{j1}$ in <N(u)> ⇒ As G is ($K_5$-e)-free, $\alpha_{j1}$ is non-adjacent to $\alpha_{i1}$, $\alpha_{i2}$ and as G is ($K_2 \cup 2K_1$)-free $\alpha_{j1}$ has no $\alpha_i$-vertex in G ⇒ j > i by **Lemma** and ∃ $\alpha_j$–vertices say $\alpha_{j2}$, $\alpha_{j3}$ ∈ N(u) not common to both a and b s.t. a$\alpha_{j2}$, b$\alpha_{j3}$ ∈ E(G) (else color all common $\alpha_j$–vertices of a, b by $\alpha_i$. If u has no $\alpha_j$–vertex, then color u by $\alpha_j$ and if u has a unique $\alpha_j$–vertex in <N(u)>, then u has |R|+1 unique colors in this coloring). As $\alpha_{j1}\alpha_{i1}$, $\alpha_{j1}\alpha_{i2}$ ∉ E(G) and G is ($K_2 \cup 2K_1$)-free, by **Lemma** $\alpha_{i1}$ ($\alpha_{i2}$) is adjacent to both $\alpha_{j2}$, $\alpha_{j3}$. Thus if $\alpha_{m1}$ is a $\alpha_m$–vertex of a, non-adjacent to b, for m ∈ {1, 2,..., χ-|R|-1}, then <a∪$\alpha_{11}$∪..∪$\alpha_{χ-|R|-11}$> = $K_{χ-|R|+1}$ ⇒ |R| > 2 (else <u∪a∪$\alpha_{11}$∪..∪$\alpha_{χ-|R|-1}$> = $K_{χ-|R|+1}$ and χ ≤ ω+1). Let x ∈ R-$R_1$. W.l.g. let x$\alpha_{i1}$ ∈ E(G) ⇒ xa ∉ E(G) (else as G is ($K_5$-e)-free x$\alpha_{m1}$ ∈ E(G) ∀ m and <u∪a∪x∪$\alpha_{11}$∪..∪$\alpha_{χ-|R|-1}$> = $K_{χ-|R|+2}$ and χ ≤ ω+1) ⇒ x$\alpha_{j2}$ ∉ E(G) ⇒ As G is ($K_2 \cup 2K_1$)-free, x$\alpha_{j1}$, x$\alpha_{j3}$ ∈ E(G) ⇒ As before x$\alpha_{i2}$, xb ∉ E(G). Let ax' ∈ E(G) ⇒ As G is ($K_2 \cup 2K_1$)-free, bx' ∈ E(G). Now x'$\alpha_{i2}$ ∈ E(G) (else <a, x', x, $\alpha_{i2}$> = ($K_2 \cup 2K_1$)) ⇒ x'$\alpha_{j3}$ ∉ E(G) (else <x', $\alpha_{j3}$, b, $\alpha_{i2}$, u> = ($K_5$-e)) ⇒ x'$\alpha_{j1}$ ∈ E(G) (else x' has no $\alpha_j$-vertex. Then color x' by $\alpha_j$; a by color of x; u by color of a). But then <x', a, b, $\alpha_{j1}$, u> = ($K_5$-e), a contradiction. Hence the **Claim** holds.

**Case B.1**: χ-|R|-1 = 1.
⇒ |R| = 3 = ω (else χ ≤ ω+1). Let x ∈ R-$R_1$. As $N_1$ = 2, by **Lemma** w.l.g. let x$\alpha_{11}$ ∈ E(G) ⇒ xa ∉ E(G) (else <u, a, x, $\alpha_{11}$> = $K_4$). Let ax' ∈ E(G) where x', x have same color.
**Case B.1.1**: $\alpha_{11}$ is a unique $\alpha_1$-vertex of a.
⇒ G has only two $\alpha_1$-vertices, $\alpha_{12}$ is a unique $\alpha_1$-vertex of b and ∃ a', b' ∈ V(G)-N(u) s.t. b'$\alpha_{11}$, a'$\alpha_{12}$ ∈ E(G) ⇒ xa', x'a ∈ E(G) and a'b' ∈ E(G) (else color b by $\alpha_1$; $\alpha_{12}$ by color of a; a', u by color of b). Also x$\alpha_{12}$ ∈ E(G) (else $\alpha_{11}$ is a unique $\alpha_1$-vertex of a and x ⇒ $\alpha_{11}$a', $\alpha_{11}$x' ∈ E(G) ⇒ x'b' ∉ E(G) {else <x', b', a', $\alpha_{11}$> = $K_4$} ⇒ b'x ∈ E(G) and <x, b', a', $\alpha_{11}$> = $K_4$) ⇒ xb ∉ E(G), bx' ∈ E(G) ⇒ a'$\alpha_{11}$, b'$\alpha_{12}$ ∉ E(G) ⇒ ab' (ba') ∈ E(G) (else let a = u', color u by color of a, then u' has |R|+1 unique colors) and x' has no $\alpha_1$–vertex in <N(u)> and hence no $\alpha_1$-vertex in G. Color x' by $\alpha_1$; a by color of x; u by color of a, a contradiction.
**Case B.1.2**: ∃ $\alpha_1$' ∈ V(G)-N(u) s.t. a$\alpha_1$' ∈ E(G).
As G is ($K_2 \cup 2K_1$)-free, b$\alpha_1$' ∈ E(G). Also x$\alpha_1$' ∈ E(G) (else x$\alpha_{12}$ ∈ E(G) ⇒ xb ∉ E(G) ⇒ bx' ∈ E(G) ⇒ x'$\alpha_1$' ∈ E(G) {else <u, x, x', $\alpha_1$'> = $K_2 \cup 2K_1$} and <a, b, x', $\alpha_1$'> = $K_4$). Let xa' ∈ E(G) where a, a' have same color. Then $\alpha_{12}$ is adjacent to either x or a' (else <x, a', a, $\alpha_{12}$> = ($K_2 \cup 2K_1$)).

**B.1.2.1**: $\alpha_{12}x$ is $\in E(G)$

$\Rightarrow xb \notin E(G)$. Let $xb' \in E(G)$ where $b'$, $b$ have same color. Now $bx'$, $x'b$, $x'a' \in E(G) \Rightarrow x'\alpha_1' \notin E(G) \Rightarrow x'\alpha_{11}$, $x'\alpha_{12} \in E(G)$. Also $\alpha_1'a' \in E(G)$ (else $<x'$, $a'$, $u$, $\alpha_1'> = (K_2 \cup 2K_1)$). Similarly, $\alpha_1'b' \in E(G) \Rightarrow a'b' \notin E(G)$ (else $<b'$, $x$, $a'$, $\alpha_1'> = K_4$) $\Rightarrow ba'$, $ab' \in E(G)$ and $a'\alpha_{12} \notin E(G)$ (else $<a'$, $\alpha_{12}$, $x'$, $b> = K_4$). Similarly, $b'\alpha_{11} \notin E(G)$ But then color $\alpha_{12}$ by color of $a$; $\alpha_{11}$ by color of $b$; $u$ by $\alpha_1$, a contradiction.

**B.1.2.2**: $\alpha_{12}x \notin E(G)$

$\Rightarrow \alpha_{12}x' \in E(G)$ (else $<a, x', \alpha_{12}, x> = (K_2 \cup 2K_1)$) and $\alpha_{12}a' \in E(G)$ (else $<a', x, \alpha_{12}, a> = (K_2 \cup 2K_1)$). Now $\exists b' \in V(G)-N(u)$ s.t. $x'b' \in E(G)$ (else $b$ is a unique vertex of $x'$ of that color $\Rightarrow ba' \in E(G)$ {else color $b$ by color of $a$; $a$ by color of $x$; $x'$, $u$ by color of $b$} and $<b, a', \alpha_{12}, x'> = K_4$). Also $a'b' \in E(G)$ (else $ab'$, $ba' \in E(G)$, $a$ is a unique vertex of $b'$ of that color $\Rightarrow x'b \in E(G)$ and $<b, \alpha_{12}, x', a'> = K_4$) $\Rightarrow b'\alpha_{12} \notin E(G)$ (else $<\alpha_{12}, b', x', a'> = K_4$) $\Rightarrow b'\alpha_{11}$, $b'\alpha_1' \in E(G)$. Now $xb' \notin E(G)$ (else as $\omega = 3$, $a'\alpha_{11} \notin E(G) \Rightarrow a'\alpha_1' \in E(G)$ and $<a', b', x, \alpha_1'> = K_4$) $\Rightarrow ab' \notin E(G)$ (else $<a, b', x, \alpha_{12}> = (K_2 \cup 2K_1)$) and $b$ is a unique vertex of $a$, $x$ of that color $\Rightarrow ba' \in E(G)$ (else color $a$, $x$ by color of $b$; $b$ by color of $a$; $u$ by color of $x$). Similarly, $bx' \in E(G)$ and $<a', b, x', \alpha_{12}> = K_4$, a contradiction.

**Case B.2**: $\chi-|R|-1 > 1$.

**Case B.2.1**: $\alpha_{11}\alpha_{21} \in E(G)$.

Now if $a\alpha_{31} \in E(G)$, then $\alpha_{31}$ is adjacent to both $\alpha_{11}$, $\alpha_{21}$ (else as $G$ is $(K_5-e)$-free, $\alpha_{31}$ is non-adjacent to both and $<\alpha_{11}, \alpha_{21}, \alpha_{31}, b> = (K_2 \cup 2K_1)$). Thus $<u \cup a \cup \{\alpha_{11}, .., \alpha_{\chi-|R|-11}\}> = K_{\chi-|R|+1} \Rightarrow \chi \leq \omega+|R|-1$ and $|R| = 3$. Let $x \in R-R_1$. Now $x$ is not adjacent to some $\alpha_{i1}$ adjacent to $a$ (else as $G$ is $(K_5-e)$-free, $<u \cup a \cup x \cup \{\alpha_{11}, .., \alpha_{\chi-|R|-11}\}> = K_{\chi-|R|+2} \Rightarrow \chi \leq \omega+|R|-2 \Rightarrow \chi \leq \omega+1$). W.l.g. let $x\alpha_{21} \notin E(G) \Rightarrow$ As $N_i = 2$, $x\alpha_{22} \in E(G) \Rightarrow xb \notin E(G)$ (else let $R_1' = \{b, x\}$, then $b$, $x$ have a common $\alpha_2$-vertex in $N(u)$, which is **Case A**) $\Rightarrow x\alpha_{11} \in E(G)$ (else $<\alpha_{11}, \alpha_{21}, b, x> = K_2 \cup 2K_1$) and $xa \notin E(G)$. Let $\{a, x', a', x\}$ be an $a$-$x$ bicolor path. Then $\{b, x', b', x\}$ is a $b$-$x$ bicolor path. Also $\chi-|R|-1 = 2$ (else if $a\alpha_{31} \in E(G)$, then as above $\alpha_{31}\alpha_{11}$, $\alpha_{31}\alpha_{21} \in E(G)$, $x\alpha_{31} \notin E(G)$ and $<\alpha_{31}, \alpha_{21}, x, b> = (K_2 \cup 2K_1)$) $\Rightarrow \omega \leq 4$ (else $\chi \leq \omega+1$).

**Case B.2.1.1**: $\alpha_{12}\alpha_{22} \in E(G)$.

Then $x\alpha_{12} \notin E(G) \Rightarrow x'\alpha_{12} \in E(G)$ (else $<a, x', x, \alpha_{12}> = (K_2 \cup 2K_1)$). Similarly, $x'\alpha_{21} \in E(G) \Rightarrow x'\alpha_{11} (x'\alpha_{22}) \notin E(G)$ (else $<u, a, \alpha_{11}, \alpha_{21}, x'> = K_5-e$) $\Rightarrow \exists \alpha_1$-vertex $\alpha_1' \in V(G)-N(u)$ s.t. $x'\alpha_1' \in E(G)$ (else color $a$, $\alpha_{12}$ by color of $x$; $x'$, $u$ by color of $a$). Similarly, $\exists \alpha_2$-vertex $\alpha_2' \in V(G)-N(u)$ s.t. $x'\alpha_2' \in E(G) \Rightarrow$ Each of $a$, $b$ is adjacent to $\alpha_i'$, $i = 1, 2$ and $<a, b, x', \alpha_1', \alpha_2'> = K_5-e$ or $K_5$, a contradiction.

**Case B.2.1.2**: $\alpha_{12}\alpha_{22} \notin E(G)$.

$\Rightarrow \alpha_{11}\alpha_{22} \in E(G)$ (else $<a, \alpha_{11}, \alpha_{22}, \alpha_{12}> = (K_2 \cup 2K_1)$). Similarly, $\alpha_{12}\alpha_{21} \in E(G) \Rightarrow x'\alpha_{21} \in E(G)$ (else $<b, x', x, \alpha_{21}> = (K_2 \cup 2K_1)$) $\Rightarrow x'\alpha_{11} \notin E(G) \Rightarrow x'\alpha_{12} \in E(G)$. Also $x\alpha_{12} \in E(G)$ (else $<x, \alpha_{22}, a, \alpha_{12}> = (K_2 \cup 2K_1)$), $b'\alpha_{21} \in E(G)$ (else $<x, b', \alpha_{21}, b> = (K_2 \cup 2K_1)$). Now $G$ has only two $\alpha_1$-vertices (else $\exists \alpha_1$-vertex say $\alpha_1' \in V(G)-N(u)$, $a\alpha_1'$, $b\alpha_1'$, $x'\alpha_1'$, $\alpha_1'\alpha_{22} \in E(G) \Rightarrow x'\alpha_{22} \notin E(G)$ {else $<a, b, x', \alpha_1', \alpha_{22}> = K_5-e$} and $\Rightarrow \exists \alpha_2$-vertex say $\alpha_2' \in V(G)-N(u)$ {else color $b$, $\alpha_{21}$ by color of $x$; $x'$ by $\alpha_2$; $u$ by color of $b$}. But then $a\alpha_2'$, $b\alpha_2'$, $x'\alpha_2' \in E(G)$ and $<a, b, x', \alpha_1', \alpha_2'> = K_5-e$ or $K_5$). Thus $\alpha_{1i}$ ($i = 1, 2$) is a unique $\alpha_1$-vertex of $a$ and $b$ resply. Then as by our assumption, $a$ has at most $|R|$ vertices of unique colors in $N(a)$ and clearly $ab' \in E(G) \Rightarrow b'\alpha_{11} \notin E(G)$ (else $<u, \alpha_{21}, a, \alpha_{11}, b'> = K_5-e$). Then color $\alpha_{11}$ by color of $b$; $a$ by $\alpha_1$; $u$ by color of $a$, a contradiction.

**Case B.2.2**: $\alpha_{11}\alpha_{21}$, $\alpha_{11}\alpha_{21} \notin E(G)$.

$\Rightarrow \alpha_{11}\alpha_{22}$, $\alpha_{12}\alpha_{21} \in E(G)$ (else $<b, \alpha_{12}, \alpha_{11}, \alpha_{21}> = (K_2 \cup 2K_1)$). Now $a$ ($b$) has at most one unique color $\alpha_i$ (else let $a$ have two vertices $\alpha_{i1}$ and $\alpha_{j1}$ with unique colors $\Rightarrow$ as $G$ is $(K_2 \cup 2K_1)$-free, $G$ has only two $\alpha_i$-vertices and $\alpha_j$-vertices. Then color $a$, $\alpha_{i2}$ by $\alpha_j$; $\alpha_{j1}$ by $\alpha_i$; $u$ by color of $a$).

**Case B.2.2.1**: $a$ has a unique $\alpha_1$-vertex $\alpha_{11} \Rightarrow b$ has a unique $\alpha_1$-vertex $\alpha_{12}$.

Let $\alpha_2' \in V(G)-N(u)$ s.t. $a\alpha_2' \in E(G) \Rightarrow b\alpha_2' \in E(G)$. Also $\alpha_2'\alpha_{11}$ ($\alpha_2'\alpha_{12}$) $\in E(G)$ (else $<\alpha_2', b, \alpha_{11}, \alpha_{21}> = (K_2 \cup 2K_1)$). Again $\exists a'$, $b'$ s.t. $b'\alpha_{11}$, $a'\alpha_{12} \in E(G) \Rightarrow b'\alpha_{21} \in E(G)$ (else $<\alpha_{11}, b', b, \alpha_{21}> = (K_2 \cup 2K_1)$). Similarly, $a'\alpha_{22} \in E(G) \Rightarrow a'b' \in E(G)$ (else color $b$ by $\alpha_1$; $\alpha_{12}$ by color of $a$; $a'$, $u$ by color of $b$) $\Rightarrow$ either $\alpha_2'a'$ or $\alpha_2'b' \in E(G)$ (else $<a', b', u, \alpha_2'> = (K_2 \cup 2K_1)$). W.l.g. let $a'\alpha_2' \in E(G)$.

**Claim 1**: $\chi-|R|-1 = 2$.

Else $\exists\ \alpha_3' \in V(G)-N(u)$ s.t. $a\alpha_3'$, $b\alpha_3'$, $\alpha_3'\alpha_{11} \in E(G) \Rightarrow \alpha_3'\alpha_2' \notin E(G)$ (else $<\alpha_{11}, a, \alpha_3', \alpha_2', b> =$ $K_5$-e). Also $a'\alpha_3' \in E(G)$ (else $<a', \alpha_2', \alpha_3', u> = (K_2 \cup 2K_1)) \Rightarrow a'\alpha_{11} \in E(G)$ (else $\alpha_{12}$ is a unique $\alpha_1$-vertex of a' and $ba' \in E(G)$ {else color b, a' by $\alpha_1$; $\alpha_{12}$ by color of a; u by color of b} $\Rightarrow <\alpha_3', b, \alpha_{12}, a', \alpha_2'> = K_5$-e) $\Rightarrow b'\alpha_2' \notin E(G)$ (else $b'\alpha_3' \in E(G)$ and $<\alpha_3', b', \alpha_{11}, a', \alpha_2'> = K_5$-e) $\Rightarrow b'\alpha_3' \notin$ E(G) and $b'\alpha_{22}$, $b'\alpha_{32} \in E(G) \Rightarrow <\alpha_{22}, b', \alpha_{11}, a', \alpha_{32}> = K_5$-e, a contradiction. Hence **Claim 1** holds.

Thus $\chi = |R|+3 \Rightarrow |R| = 3$ (else $\omega = 3$ {else $\chi \leq \omega+1$} and $ba' \notin E(G) \Rightarrow$ a is a unique vertex of b of that color $\Rightarrow \alpha_{11}a' \in E(G)$ {else color b, $\alpha_{11}$ by color of a; a by $\alpha_1$; u by color of b} $\Rightarrow b'\alpha_2' \notin E(G)$ and $b'\alpha_{22} \in E(G) \Rightarrow <\alpha_{22}, b', \alpha_{11}, a'> = K_4$) and $\omega \leq 4$. Let $x \in R-R_1$ and w.l.g. $xa \notin E(G)$. Let {a, x', a', x} be a a-x bicolor path. Then $ab' \in E(G)$ (else color u by color of a, let a = u', then u' has $|R|+1$ unique colors in this coloring).

**Claim 2**: $xb \notin E(G)$.

Let if possible $xb \in E(G)$. Then $x\alpha_{12} \notin E(G)$ (else if $R_1' = \{b, x\}$, then this is **Case A**) and $\alpha_{11}$ is a unique $\alpha_1$-vertex of a, x $\Rightarrow \alpha_{11}a'$, $\alpha_{11}x' \in E(G)$ (else color a, x by $\alpha_1$; $\alpha_{11}$ by color of a (or x). Also $\alpha_{12}x' \in E(G)$ (else color $\alpha_{12}$ by color of x; b by $\alpha_1$; u by color of b). Again $x'b' \in E(G)$ (else b is a unique vertex of x' of that color $\Rightarrow ba' \in E(G) \Rightarrow \alpha_2'x', \alpha_2'a' \notin E(G)$ {$<\alpha_2', \alpha_{12}, b, a', x'> = K_5$ or $K_5$-e} and $<x', a', u, \alpha_2'> = (K_2 \cup 2K_1))$. Now as G is $(K_2 \cup 2K_1)$-free $\alpha_2'$ is adjacent to at least two of $\{a', b', x'\} \Rightarrow <\alpha_2', \alpha_{11}, b', a', x'> = K_5$ or $K_5$-e, a contradiction. Hence **Claim 2** holds.

Thus b has two vertices of color of a and $ba' \in E(G)$. Also {b, x', b', x} is a b-x bicolor path. Now $x'\alpha_2', x'\alpha_{12} \notin E(G)$ (else $<x', b, a', \alpha_2', \alpha_{12}> = K_5$ or $K_5$-e) $\Rightarrow x\alpha_2' \in E(G)$ (else $<u, x, x', \alpha_2'> = (K_2 \cup 2K_1)$, $x\alpha_{12} \in E(G)$ (else color $\alpha_{12}$ by color of x; b by $\alpha_1$; u by color of b) and $\alpha_2'b' \in E(G)$ (else $<x', b', u, \alpha_2'> = (K_2 \cup 2K_1)$). But then $<b', \alpha_2', a', x, \alpha_{12}> = K_5$ or $K_5$-e, a contradiction.

**Case B.2.2.2**: a has $\alpha_i$-vertex $\alpha_i' \in V(G)-N(u)$, $\forall\ i = 1,..., \chi-|R|-1$.

$\Rightarrow$ As before $b\alpha_i' \in E(G)$, $\alpha_i'\alpha_{jk} \in E(G)$ for $j \neq i$, $k = 1, 2$. Again $\alpha_{11}$ or $\alpha_{21}$ has a vertex b' of color of b (else color $\alpha_{11}$ or $\alpha_{21}$ by color of b; $\alpha_{12}$ by $\alpha_2$; $\alpha_2'$, u by $\alpha_1$) and $\alpha_{11}b' \in E(G)$ iff $\alpha_{21}b' \in E(G)$ (else $<\alpha_{11}, b', b, \alpha_{21}> = K_2 \cup 2K_1) \Rightarrow \alpha_{11}b', \alpha_{21}b', \alpha_{12}a', \alpha_{22}a' \in E(G)$ where a', a have same color. Let {a, x', a', x} be a bicolor a-x path.

**Claim**: $\alpha_1'\alpha_2' \notin E(G)$.

Else x' is non-adjacent to $\alpha_1'$ or $\alpha_2'$ {else $<b, a, \alpha_1', \alpha_2', x'> = K_5$ or $K_5$-e} and adjacent to the other {else $<\alpha_1', \alpha_2', x', u> = K_2 \cup 2K_1$}. Let $x'\alpha_1' \in E(G)$ and $x'\alpha_2' \notin E(G) \Rightarrow x'b \notin E(G)$ {else $<\alpha_2', \alpha_1', a, b, x'> K_5$-e} and $bx \in E(G)$, $x\alpha_1' \notin E(G)$ and $x\alpha_{12} \in E(G) \Rightarrow <b, x, \alpha_2', \alpha_{12}, a'> = K_5$ or $K_5$-e, a contradiction. Hence the **Claim** holds.

a) $x'\alpha_1' \notin E(G) \Rightarrow x'\alpha_2' \notin E(G)$.

$\Rightarrow x'\alpha_{11}, x'\alpha_{12} \in E(G)$ and $x'b \in E(G)$ (else b is a unique vertex of x' of that color and $ba' \in E(G)$ {else color b by color of a; a by color of x; x', u by color of b}. Also $\alpha_2'a' \in E(G)$ {else $<x', a', \alpha_2', u> = K_2 \cup 2K_1$} $\Rightarrow <\alpha_2', b, \alpha_{12}, a', x'> = K_5$ or $K_5$-e) $\Rightarrow ab' \notin E(G)$ (else $<x', a, \alpha_{11}, b', \alpha_2'> = K_5$-e), $a'b' \in E(G)$ (else $<u, a, a', b'> = K_2 \cup 2K_1$) and b is a unique vertex of a of that color. Also $x\alpha_1', x\alpha_2' \in E(G)$ (else $<\alpha_1', \alpha_2', x, u> = K_2 \cup 2K_1) \Rightarrow xb' \notin E(G)$ (else $<\alpha_1', x, a', b', \alpha_2'> = K_5$-e). Thus b is a unique vertex of a, x of that color $\Rightarrow bx' \in E(G)$ (else color a, x by color of b; b by color of x; u by color of a). Similarly, $ba' \in E(G)$ and $<\alpha_2', b, a', \alpha_{12}, x'> = K_5$-e, a contradiction.

b) $x'\alpha_1' \in E(G) \Rightarrow x'\alpha_2' \in E(G)$.

$\Rightarrow x'b \notin E(G)$ (else $<\alpha_1', x', a, b, \alpha_2'> = K_5$-e) $\Rightarrow bx \in E(G) \Rightarrow x\alpha_{12} \notin E(G)$ (else can be proved as **Case A**) $\Rightarrow x\alpha_1' \in E(G) \Rightarrow x\alpha_2' \in E(G)$. Also $x'\alpha_{12} \in E(G)$ (else $<a, x', \alpha_{12}, x> = K_2 \cup 2K_1$) and $a'b' \in E(G)$ (else color b by color of x; x by color of a; a', u by color of b). Now $a'\alpha_2' \notin E(G)$ (else $b'\alpha_2'$, $b'\alpha_{12} \notin E(G)$ {else $<b', a', x, \alpha_{12}, \alpha_2'> = K_5$-e or $K_5$} $\Rightarrow b'\alpha_1' \in E(G)$ and $<b', \alpha_1', u, \alpha_2'> = K_2 \cup 2K_1) \Rightarrow a'\alpha_1' \notin E(G)$, $\alpha_1'b', \alpha_2'b', a'\alpha_{11} \in E(G)$ and $b'x \notin E(G)$ (else $<\alpha_2', \alpha_{11}, x, b', a'> = K_5$-

e) ⇒ ab' ∈ E(G) (else b is a unique vertex of a, x of that color. Color a, x by color of b; b by color of x; u by color of a) ⇒ <α$_1$', x', a, b', α$_2$'> = $K_5$-e, a contradiction.

This proves **Theorem 4**.

The following Table summarizes the results.

| G free of | χ ≤ |
|---|---|
| {($P_4$+$K_1$), Chair}; {($P_4$+$K_1$), $K_{1,3}$} | 2ω-1 |
| {($P_4$+$K_1$), $P_3$U$K_1$} | |
| {($P_4$+$K_1$), ($K_2$U2$K_1$)} | $\frac{3\omega}{2}$ |
| {HVN}, Chair}; {HVN, $K_{1,3}$} | |
| {HVN, $P_3$U$K_1$} | |
| {HVN, ($K_2$U2$K_1$)} | |
| {$K_4$, Chair} | ω+1 |
| {$K_5$-e, $P_3$U$K_1$} | |
| {$K_5$-e, ($K_2$U2$K_1$)} | |

**Table**

**Note:** In [7], Kierstead and Schmerl proved that if G induces neither $K_{1,3}$ nor $K_5$-e, then χ(G) ≤ ω+1.

**Examples to show that Exclusion of certain graphs is Necessary:**
1. H = $C_5$+ $C_5$+…+ $C_5$ (m times, m >1). Here ω(H) = 2m, χ(H) = 3m. H is free of each of $P_3$U$K_1$, $K_2$U2$K_1$ but has both HVN, $K_5$-e induced and χ > ω+1.
2. Grotzsch graph has ω = 2, χ = 4, has none of $P_4$+$K_1$, HVN, $K_5$-e but has a Chair induced and χ > $\frac{3\omega}{2}$.